\documentclass[11pt]{amsart}

\usepackage[T1]{fontenc}
\usepackage{lmodern}
\usepackage[margin=1in]{geometry}
\usepackage{mathtools,amssymb,array,microtype}
\usepackage[hidelinks]{hyperref}

\newtheorem{theorem}{Theorem}[section]
\newtheorem{lemma}[theorem]{Lemma}
\newtheorem{corollary}[theorem]{Corollary}
\theoremstyle{definition}
\newtheorem{example}[theorem]{Example}
\newtheorem{assumption}{Assumption}

\numberwithin{equation}{section}

\newcommand{\dd}{\mathrm d}
\newcommand{\R}{\mathbb R}
\newcommand{\Pp}{\mathbb P}
\newcommand{\E}{\mathbb E}
\newcommand{\ind}{\mathbf 1}
\newcommand{\diag}{\operatorname{diag}}
\newcommand{\op}{o_{\mathrm p}}
\newcommand{\Op}{O_{\mathrm p}}
\newcommand{\Normal}{\mathcal N}
\newcommand{\assref}[1]{\hyperref[ass:#1]{\textup{(#1)}}}

\title[Infinite dimensional limits for the giant component]
{Infinite dimensional Gaussian limits for the giant component of inhomogeneous random graphs}
\author[H. Ma]{Haoxuan Ma}
\address{Sorbonne Universit\'e, France}
\email{Haoxuan.Ma@etu.sorbonne-universite.fr}
\author[W. Sun]{Wen Sun}
\address{University of Science and Technology of China, China}
\email{wensun.ustc@gmail.com}
\date{November 6, 2026}
\hypersetup{pdftitle={Infinite dimensional Gaussian limits for the giant component of inhomogeneous random graphs},
pdfauthor={Haoxuan Ma and Wen Sun}}

\begin{document}
\begin{abstract}
We prove Gaussian limits for the type measure of the giant component
in sparse inhomogeneous random graphs. Countable types yield convergence
in a weighted Hilbert space under boundedness, irreducibility,
supercriticality, and a uniform tail condition on the empirical type
weights. For types on a finite dimensional closed Riemannian manifold,
smooth kernels yield convergence in a negative Sobolev space for any
Borel type distribution. We allow deterministic types and first order
kernel perturbations; independent type sampling adds an independent
Gaussian term. The proof combines a seeded exploration identity,
Hilbert space martingale limits, tail estimates, and inversion of the
branching operator outside the giant. Finite types with random parameters
are covered. A heat kernel example gives explicit variances for all
spectral modes.
\end{abstract}
\maketitle

\section{Introduction}\label{sec:introduction}

We consider a sparse inhomogeneous random graph $G_n$ on
$[n]=\{1,\ldots,n\}$.  Each vertex $v$ has a type
$x_v^{(n)}$ in a space $\mathcal S$.  We consider finite and countably
infinite type sets, and finite dimensional closed smooth Riemannian
manifolds.  Type sets carry the discrete topology and manifolds their
manifold topology; both carry the corresponding Borel $\sigma$-algebra.
  Given the types and a symmetric nonnegative measurable kernel
$\kappa_n:\mathcal S\times\mathcal S\to[0,\infty)$, distinct unordered
vertex pairs are joined independently with probabilities
\begin{equation}
 p_n(x_v^{(n)},x_w^{(n)})
 =1-\exp\{-\kappa_n(x_v^{(n)},x_w^{(n)})/n\},\qquad v\ne w.
 \label{eq:common-graph-model}
\end{equation}
We write $v\sim_n w$ when the vertices are joined.  There are no self-loops
or multiple edges.  The vertex types and the kernel are
deterministic in the two main theorems; random types and parameters are
treated in the corollaries below.

Let $C_1^{(n)}$ and $C_2^{(n)}$ be the largest and second largest connected
components, ordered by their numbers of vertices. We break ties by a
fixed rule and put $C_2^{(n)}=\varnothing$ if the graph is connected.
Writing $\delta_x$ for the unit point mass at $x$, define the empirical
type distribution and the type counting measure of the largest component by
\begin{equation}
 \mu^{(n)}=\frac1n\sum_{v=1}^n\delta_{x_v^{(n)}},\qquad
 \mathcal L_n=\sum_{v\in C_1^{(n)}}\delta_{x_v^{(n)}}.
 \label{eq:common-type-measures}
\end{equation}
For a bounded real valued measurable function $f$ on $\mathcal S$, we write
\[
 \mathcal L_n(f)=\int_{\mathcal S}f\,\dd\mathcal L_n
 =\sum_{v\in C_1^{(n)}}f(x_v^{(n)}).
\]
In particular, $\mathcal L_n(1)=|C_1^{(n)}|$ is the total giant size.
We study Gaussian fluctuations of the entire type measure on the
$\sqrt n$ scale.  For countably many types we view it as a vector in a
weighted Hilbert space; for manifold types we use a negative Sobolev
space.

The phase transition in the sparse Erd\H{o}s--R\'enyi graph
\cite{ErdosRenyi1960} is the one type case of this problem.  In the
supercritical regime, the largest component contains a positive fraction
of the vertices.  Its central limit theorem has been proved by several
methods; see Stepanov \cite{Stepanov1970}, Martin-L\"of
\cite{MartinLof1986}, Pittel \cite{Pittel1990}, and Barraez, Boucheron,
and Fernandez de la Vega \cite{BarraezBoucheronVega2000}.
Puhalskii \cite{Puhalskii2005} treats normal, moderate, and large
deviations through stochastic process methods, and a recent proof of the
central limit theorem by local approximation is given in
\cite{AnastosErdeKangPfenninger2026}.  Moderate deviation principles for
the giant size are also proved by Ameskamp and L\"owe
\cite{AmeskampLowe2011} and Sun \cite{Sun2023}.  Large deviation
principles are obtained by O'Connell \cite{OConnell1998} for the largest
component and by Andreis, K\"onig, and Patterson
\cite{AndreisKoenigPatterson2021} for the collection of component sizes.  Process limits for the giant
size in dynamic Erd\H{o}s--R\'enyi graphs are proved by Enriquez, Faraud,
and Lemaire \cite{EnriquezFaraudLemaire2025} and by Corujo, Lemaire,
and Limic \cite{CorujoLemaireLimic2026}.

The finite type kernel graphs of S\"oderberg \cite{Soderberg2002}
extend the one type model by allowing different interaction rates between
groups of vertices. For general type spaces, Bollob\'as, Janson, and Riordan
\cite[Theorems~3.1 and~9.10]{BollobasJansonRiordan2007} locate the phase
transition at branching operator norm one and prove laws of large numbers
for the size and type composition of the largest component.  Its limiting
type measure is determined by the survival probabilities of an associated
Poisson branching process.  This first order description is the starting
point of our fluctuation analysis. With suitable rank one kernels, the
framework also includes expected degree and Poissonian models
\cite{ChungLu2002,NorrosReittu2006}.  The bounded-kernel assumptions
used here cover the bounded-weight cases of these models.

For finitely many types, the centered vector central limit theorem,
including the covariance arising from independent type allocation, can
be obtained from Neal's multitype epidemic theorem
\cite[Theorems~5.5 and~6.2; Section~7.3]{Neal2006}.  Clancy
\cite{Clancy2025} includes deterministic perturbations of the type
proportions and kernel on the $n^{-1/2}$ scale, with the limit stated
after applying the kernel to the type vector.  Bhamidi, Budhiraja, and
Sakanaveeti \cite{BhamidiBudhirajaSakanaveeti2026} prove functional
central limit theorems for microscopic component densities and joint
macroscopic observables, including the giant type vector, for finitely
many vertex types.  Yu and Sun \cite{YuSun2025} establish moderate
deviation principles for the giant type vector under additional
assumptions.  Andreis, K\"onig, Langhammer, and Patterson
\cite{AndreisKoenigLanghammerPatterson2023} prove a joint large
deviation principle for component configurations, including their vertex
types; their result also covers compact metric type spaces.

For rank one Poisson graphs, Clancy \cite{ClancyPoisson2026} proves a
joint process limit for the giant size and its weighted volume under
convergence of the empirical weight distribution and its second moment.
The assumptions allow deterministic weight arrays and unbounded weights,
and the centering uses the empirical fixed point.  Our results instead
concern the entire type measure at a fixed kernel, with no rank
restriction, under boundedness and the regularity assumptions below.

Infinite dimensional Gaussian approximations also arise in spatial
epidemic models.  Andersson and Djehiche
\cite{AnderssonDjehiche1995} obtain Hilbert space valued process limits
for multitype approximations to Kendall's spatial epidemic, with an
infinite dimensional Ornstein--Uhlenbeck fluctuation limit.  They also
study the final epidemic size when both the number of types and the
number of individuals of each type tend to infinity
\cite{AnderssonDjehiche1997}.

Our two main results concern countably many types and types on a closed
manifold.  They give convergence of the entire giant type measure in a
weighted Hilbert space and a negative Sobolev space, respectively.  Each
result therefore yields joint Gaussian limits for all continuous linear
observations of that measure, including their cross-covariances. This
covers bounded type weights in the countable model and the Sobolev test
functions specified below on a manifold. In particular, size and type
composition can be studied jointly through the same Gaussian limit.
The proofs combine finite dimensional martingale limits with coordinate
or spectral tail estimates and inversion of the linearized branching
equation.  The results
include kernels of infinite rank; finite types are covered by taking a
finite index set in the countable model.
A finite type CLT alone does not supply convergence of the full type
measure: one must also control the omitted coordinates and the inverse
of the limiting derivative.  We prove these bounds directly in the
weighted Hilbert and negative Sobolev spaces.  The two inverse arguments
are instances of the same factorization through the branching operator
outside the giant.  Ratios of continuous linear observations then give
Gaussian limits for the type proportions within the giant; the
calculation is given in the
\hyperref[rem:type-proportions]{remark on joint observables and type proportions}
in Subsection~\ref{sec:compact-model}.

The boundedness assumption is essential to the estimates used in this
proof.  It gives a uniform positive lower bound on the exploration
factors, bounded fourth moments of the increments, and uniform control
of their conditional variances.  Unbounded kernels, including rank one
models with unbounded weights, require additional truncation and
integrability arguments.  On a manifold, the Sobolev framework provides
compactness of point masses, bounded multiplication, and smoothing by
the kernel.  A general compact metric space would require a specified
function space with corresponding properties; compactness of the type
space alone does not provide them.

\subsection{Notation}\label{sec:notation}

We collect the spaces and operations used throughout the paper.
The assumptions for the two type spaces are specified in
Section~\ref{sec:model}, and the exploration processes are introduced
in Section~\ref{sec:proof}.

\paragraph{Probability and convergence.}
We work over the real numbers and write $\mathbb N=\{1,2,\ldots\}$,
$\Pp$ and $\E$ for probability and expectation,
$\ind_A$ for the indicator of $A$, and $\mathrm i=\sqrt{-1}$.
The symbol $\boxtimes$ denotes the product of measures or probability laws.
All limits are taken as $n\to\infty$ unless stated otherwise.
The symbols $\Rightarrow$ and $\xrightarrow{\mathrm p}$ denote weak
convergence of Borel laws and convergence in probability in the stated
normed space.  For positive deterministic $a_n$, $X_n=\Op(a_n)$ means
that $\|X_n\|/a_n$ is bounded in probability, and $X_n=\op(a_n)$ means
that this ratio converges to zero in probability.
The notation $\Normal_d(m,\Sigma)$ denotes the Gaussian law on $\R^d$
with mean $m$ and covariance matrix $\Sigma$; $\Normal(m,\Sigma)$ denotes
the corresponding law on a real Hilbert space, with covariance operator
$\Sigma$ defined with respect to the inner product of that space.  For a discrete-time process,
$\Delta X_{n,t}=X_{n,t+1}-X_{n,t}$.
Constants denoted by $C$ may change from line to line; dependence on
fixed quantities is indicated when needed. Unless explicitly indicated,
they do not depend on $n$, the exploration time, or the type coordinate.

\paragraph{Sequences and matrices.}
The index set $\mathcal I$ is either $\mathbb N$ or $[K]$ for a fixed
positive integer $K$, and $\mathcal I_m=\{i\in\mathcal I:i\le m\}$.
Sums over types run over $\mathcal I$.
For functions, $\|f\|_\infty=\sup_x|f(x)|$.
For sequences, $\|x\|_1=\sum_i|x_i|$,
$\|x\|_2^2=\sum_i x_i^2$, and $\|x\|_\infty=\sup_i|x_i|$;
$\ell^1,\ell^2,\ell^\infty$ are the spaces on which these norms are finite.
For positive weights $w_i\ge1$, put
\[
 \mathcal H_w=\left\{x:\|x\|_w^2=\sum_iw_ix_i^2<\infty\right\},
 \qquad \langle x,y\rangle_w=\sum_iw_ix_iy_i,
 \qquad C_w=\sum_iw_i^{-1}.
\]
The coordinate unit vectors are
$(\mathbf e_j)_i=\ind_{\{i=j\}}$, so
$(\mathbf e_j/\sqrt{w_j})_j$ is an orthonormal basis of $\mathcal H_w$.
The coordinate projections are
\begin{equation}
 P_mx=(x_i\ind_{\{i\in\mathcal I_m\}})_i.
 \label{eq:countable-projection}
\end{equation}
Products, powers, and exponentials of sequences act coordinatewise;
$(a\odot b)_i=a_ib_i$ and $(e^a)_i=e^{a_i}$.
We also write $ab=a\odot b$, and $D_a=\diag(a)$.
For a matrix or countable array $B$, $(Bx)_i=\sum_jB_{ij}x_j$
whenever the series converges absolutely.  Its weighted row norm is
$\|B\|_{w,\infty}=\sup_i(\sum_jB_{ij}^2/w_j)^{1/2}$.
The space $\operatorname{Sym}_K$ consists of symmetric real
$K\times K$ matrices and carries the Frobenius norm
$\|B\|_{\mathrm F}=(\sum_{i,j}B_{ij}^2)^{1/2}$.
Matrix transpose is denoted by $\mathsf T$, with
$B^{-\mathsf T}=(B^{-1})^{\mathsf T}$ for invertible $B$.
On $\R^K$ the norm $\|\cdot\|_2$ is Euclidean, and
$\operatorname{tr}(B^{\mathsf T}A)=\sum_{i,j}B_{ij}A_{ij}$.

\paragraph{Linear operators.}
For a bounded linear map $A:E\to F$ between normed spaces,
$\|A\|_{\mathrm{op}}=\sup_{\|x\|_E\le1}\|Ax\|_F$; the domain and range
are displayed when needed.  We also write $\|A\|$ when they are clear,
and $I$ for the identity.  The notation $Dg$ denotes the Fr\'echet
derivative of a map $g$.
For a bounded linear map $A$ between real Hilbert spaces,
$A^\dagger$ denotes its Hilbert adjoint.  On a real Hilbert space
$\mathcal H$, with orthonormal basis $(u_j)$,
$\|A\|_{\mathrm{HS}}^2=\sum_j\|Au_j\|_{\mathcal H}^2$ is the
Hilbert--Schmidt norm, and $\|A\|_{\mathrm{tr}}$ is the trace norm.
In particular,
$\|A\|_{\mathrm{HS}}^2=\sum_j\|A\mathbf e_j/\sqrt{w_j}\|_w^2$
on $\mathcal H_w$.
For $u,v\in\mathcal H$,
$(u\otimes v)x=\langle x,v\rangle_{\mathcal H}u$.
On $\mathcal H_w$ we write $\otimes_w$ and $A^{*_{w}}$ for this
rank one operation and the Hilbert adjoint:
$\langle Ax,y\rangle_w=\langle x,A^{*_{w}}y\rangle_w$.

\paragraph{Measures and kernels.}
For a finite signed measure $\eta$, $\eta(f)=\int f\,\dd\eta$,
$\|\eta\|_{\mathrm{TV}}=\sup_{\|f\|_\infty\le1}|\eta(f)|$ is its total
variation norm, with the supremum over bounded measurable $f$.  Multiplication by a function is defined by
$(a\eta)(f)=\eta(af)$; the same formula defines multiplication of
a distribution when $f\mapsto af$ preserves its test-function space.
A measure on $\mathcal I$ is identified with its mass vector, so kernel
action on such a measure agrees with matrix multiplication.
For a positive measure $\mu$, $L^2(\mu)$ has norm
$\|f\|_{L^2(\mu)}^2=\int f^2\,\dd\mu$; on $\mathcal I$ we write
$\ell^2(\mu)$ for the same space.
For a kernel $k$, write
\[
 (\mathsf K_k\eta)(x)=\eta(k(x,\cdot)),\qquad
 (T_kf)(x)=\int k(x,y)f(y)\,\mu(\dd y),
 \qquad \mathsf K=\mathsf K_\kappa,
\]
whenever the pairings or integrals exist.  For finite measures,
$(\mathsf K_k\eta)(x)=\int k(x,y)\,\eta(\dd y)$.

\subsection*{Outline}

Section~\ref{sec:model} specifies the two type settings and states the
Gaussian limits and examples. Section~\ref{sec:proof} develops the seeded
exploration, its martingales, the mass balance equation, and local
linearization. Sections~\ref{sec:countable-proof} and~\ref{sec:compact}
verify the required estimates for countable and manifold types,
respectively.

\section{Model and main results}\label{sec:model}

We use the graph and giant type measure from
\eqref{eq:common-graph-model}--\eqref{eq:common-type-measures}.
Let $\mathcal E_n$ be the sigma algebra generated by the vertex types
and the kernel; it is trivial when both are deterministic.  Edges are
sampled independently conditional on $\mathcal E_n$ with the
probabilities in \eqref{eq:common-graph-model}.
The limiting type law and kernel are denoted by $\mu$ and $\kappa$.
We write $\Delta_n=\kappa_n-\kappa$ and
\begin{equation}
 H_n=\sqrt n(\mu^{(n)}-\mu)
 \label{eq:countable-iid-empirical-fluctuation}
\end{equation}
for the empirical type fluctuation.  Whenever uniform boundedness of a
kernel family is assumed, $\kappa_*<\infty$ denotes a common
deterministic upper bound for that family and its limit.  The two
subsections below specify the assumptions and the space in which the
giant type measure converges.

\subsection{Countable type spaces}\label{sec:countable}

We first describe the fluctuations of all type counts at once.
Take $\mathcal S=\mathcal I$, with the discrete spaces and operations
of Subsection~\ref{sec:notation}.  Fix a probability vector
$\mu=(\mu_i)_i$ with $\mu_i>0$ and weights $w_i\geq1$.
For the graph in \eqref{eq:common-graph-model}, write
\[
 \mu_i^{(n)}=\mu^{(n)}(\{i\}),\qquad
 L_{n,i}=\mathcal L_n(\{i\}),\qquad
 \kappa_{n,ij}=\kappa_n(i,j).
\]
We identify $\mu^{(n)}$ with its mass vector $(\mu_i^{(n)})_i$ and let $L_n=(L_{n,i})_i$.
Thus $n_{i}^{(n)}=n\mu_i^{(n)}$ is the number of type $i$ vertices, and $L_n$
is the type count vector of the largest component.
Fix a bounded symmetric nonnegative kernel $\kappa$ that is irreducible:
no $A\subset\mathcal I$ with $0<\sum_{i\in A}\mu_i<1$ has
$\kappa_{ij}=0$ for every $i\in A,j\notin A$.
Assume that it is supercritical, meaning
$\|T_\kappa\|_{\ell^2(\mu)\to\ell^2(\mu)}>1$.  A type $i$ particle in the associated branching
process has independent Poisson numbers of type $j$ children with means
$\kappa_{ij}\mu_j$.  Its survival probabilities form the maximal solution
$\rho\in[0,1]^{\mathcal I}$ of
\begin{equation}
 \rho_i=1-\exp\left\{-\sum_j\kappa_{ij}\mu_j\rho_j\right\}.
 \label{eq:countable-fixed-point}
\end{equation}

Write $\bar\rho=\sum_i\mu_i\rho_i$ for the limiting fraction of vertices
in the giant.

Under the parameter assumptions of the theorem below, the laws of large
numbers of Bollob\'as, Janson, and Riordan
\cite[Theorems~3.1 and~9.10]{BollobasJansonRiordan2007} give
\[
 \frac{L_{n,i}}n\xrightarrow{\mathrm p}\mu_i\rho_i\quad(i\in\mathcal I),
 \qquad
 \frac{|C_1^{(n)}|}n\xrightarrow{\mathrm p}\sum_i \mu_i\rho_i>0.
\]
The next theorem identifies the fluctuations around $\mu\rho$.

\begin{theorem}[Gaussian limit for the type vector]
\label{thm:countable-concrete}
Suppose that the types and kernels are deterministic, that
\begin{equation}
 \sum_iw_i\mu_i<\infty,\qquad \sum_iw_i^{-1}<\infty,
 \label{eq:main-countable-summability}
\end{equation}
and that, for some $\beta\in\mathcal H_w$,
\begin{equation}
 \sqrt n(\mu^{(n)}-\mu)\longrightarrow\beta\quad\hbox{in }\mathcal H_w,
 \qquad
 \lim_{m\to\infty}\sup_n\sum_{i>m}w_i\mu_i^{(n)}=0.
 \label{eq:main-countable-counts}
\end{equation}
Assume that the kernels are uniformly bounded and
\begin{equation}
 \kappa_n=\kappa+n^{-1/2}\Lambda+E^{(n)},\qquad
 \|\Lambda\|_{w,\infty}<\infty,\qquad
 \sqrt n\,\|E^{(n)}\|_{w,\infty}\longrightarrow0,
 \label{eq:main-countable-kernels}
\end{equation}
where $\Lambda$ is symmetric.  Then $J=I-D_{\mu(1-\rho)}\kappa$ is boundedly invertible on
$\mathcal H_w$,
and, in $\mathcal H_w$,
\begin{equation}
 \frac{L_n-n\mu\rho}{\sqrt n}
 \Rightarrow\Normal\!\bigl(
 J^{-1}\{D_\rho\beta+D_{\mu(1-\rho)}\Lambda(\mu\rho)\},
 J^{-1}D_{w\mu\rho(1-\rho)}(J^{-1})^{*_{w}}\bigr).
 \label{eq:concrete-countable-limit}
\end{equation}
\end{theorem}

The parameter perturbations contribute the first order shift
$D_\rho\beta+D_{\mu(1-\rho)}\Lambda(\mu\rho)$, while edge randomness
contributes covariance $D_{w\mu\rho(1-\rho)}$.
Both are transformed by $J^{-1}$, since $J=Dh(\mu\rho)$ is the
derivative of the mass balance map in \eqref{eq:exploration-residuals}.

For $\mathcal I=[K]$ and $w_i=1$, the summability and tail conditions
are automatic; the same theorem includes the finite type vector CLT.

When types are sampled independently, their empirical proportions also
fluctuate on the $n^{-1/2}$ scale.  The next corollary adds this source
of randomness to the edge noise.

\begin{corollary}[Independent type sampling]
\label{cor:countable-iid}
Let $\mathcal I$, $\mu$, and $w$ be as in
Subsection~\ref{sec:countable}.  Suppose that the bounded
symmetric nonnegative kernel $\kappa$ is irreducible and supercritical,
and that \eqref{eq:main-countable-summability} holds.  Let
$X_1,X_2,\ldots$ be independent with law $\mu$, and construct the graph with types
$X_1,\ldots,X_n$ and the fixed kernel $\kappa$.  Then, in $\mathcal H_w$,
\begin{equation}
 \frac{L_n-n\mu\rho}{\sqrt n}
 \Rightarrow\Normal\!\bigl(0,
 J^{-1}\{D_{w\mu\rho}-(\mu\rho)\otimes_w(\mu\rho)\}
 (J^{-1})^{*_{w}}\bigr).
 \label{eq:countable-iid-limit}
\end{equation}
In particular, for $\mathcal I=[K]$ and $w_i=1$,
\begin{equation}
 \frac{L_n-n\mu\rho}{\sqrt n}
 \Rightarrow\Normal_K\!\left(0,
   J^{-1}(D_{\mu\rho}-(\mu\rho)(\mu\rho)^{\mathsf T})J^{-\mathsf T}\right).
 \label{eq:finite-iid-limit}
\end{equation}

\end{corollary}

\begin{corollary}[Finite types with random parameters]
\label{cor:finite-random}
Take $\mathcal I=[K]$ in the discrete model, with the
limiting type probabilities and kernel satisfying its irreducibility and
supercriticality assumptions.  The partition and the symmetric nonnegative kernel
$\kappa_n$ may now be random.  The edges are sampled conditionally on $\mathcal E_n$ as in
\eqref{eq:common-graph-model}.
With $H_n$ from \eqref{eq:countable-iid-empirical-fluctuation}, put
\begin{equation}
 A_n=\sqrt n(\kappa_n-\kappa).
 \label{eq:finite-parameter-fluctuations}
\end{equation}
If
\[
 (H_n,A_n)\Rightarrow(H,A)
 \quad\hbox{in }\mathbb R^K\times\operatorname{Sym}_K,
\]
then $(L_n-n\mu\rho)/\sqrt n\Rightarrow Y$, where
\begin{equation}
 Y\mid(H,A)\sim\Normal_K\!\bigl(
 J^{-1}\{D_\rho H+D_{\mu(1-\rho)}A(\mu\rho)\},
 J^{-1}D_{\mu\rho(1-\rho)}J^{-\mathsf T}\bigr).
 \label{eq:finite-random-limit}
\end{equation}
The limiting law is the mixture of these Gaussian distributions over the
joint law of $(H,A)$.

\end{corollary}

\begin{example}[Countable types and infinite rank kernels]
\label{ex:countable-bounded-perturbations}
Choose either of the following type distributions, with the indicated
weights and cutoffs:
\[
\begin{array}{c|ccc}
 &\mu_i&w_i&m_n\\ \hline
 \text{geometric}
 &2^{-i}&2^{\alpha i}&\max\{1,\lfloor\log_2n\rfloor\}\\
 \text{polynomial}
 &Z_p^{-1}i^{-p}&i^a&\lfloor n^{1/p}\rfloor
\end{array}
\]
Here $0<\alpha<1$, $p>2$, $1<a<p-1$, and
$Z_p=\sum_{j\ge1}j^{-p}$.  In both cases, set
\[
 n_i^{(n)}=\lfloor n\mu_i\rfloor\ (2\le i\le m_n),\qquad
 n_1^{(n)}=n-\sum_{i=2}^{m_n}n_i^{(n)},\qquad
 n_i^{(n)}=0\ (i>m_n).
\]
For $\lambda>1$, $\varepsilon\ge0$, and a symmetric array
$0\le r_{ij}\le1$, take
$\kappa_{n,ij}=\kappa_{ij}=\lambda+\varepsilon r_{ij}$.
Theorem~\ref{thm:countable-concrete} applies with $\beta=\Lambda=0$ and gives
\[
 \frac{L_n-n\mu\rho}{\sqrt n}
 \Rightarrow\Normal\!\left(0,
 J^{-1}D_{w\mu\rho(1-\rho)}(J^{-1})^{*_{w}}\right)
 \quad\hbox{in }\mathcal H_w.
\]
For $\varepsilon>0$, the choice $r_{ij}=\ind_{\{i=j\}}$ has infinite rank.
\end{example}

\subsection{Types on closed manifolds}\label{sec:compact-model}

Let $\mathcal S=M$, where $(M,g)$ is a closed smooth Riemannian manifold
of dimension $d\ge0$; closed means compact and without boundary.
The manifold may be disconnected.  Let $\mu$ be any Borel probability
measure on $M$, including an atomic or singular measure, and normalize
Riemannian volume $\mathrm{vol}_g$ to mass one.  We use a negative Sobolev
space so that atomic type measures and their fluctuations belong to the
same Hilbert space.

Fix an integer $s>d/2$ and a real orthonormal Laplace eigenbasis $(e_j)$
in $L^2(\mathrm{vol}_g)$, with $-\Delta_ge_j=\lambda_je_j$, eigenvalues
in nondecreasing order with multiplicities, and $e_0=1$.
The basis and the sums below are finite for $d=0$.
Let $H^s=H^s(M)$ be the completion of $C^\infty(M)$ in the norm
below, and let $H^{-s}=H^{-s}(M)$ be its continuous dual:
\[
 \|f\|_s^2=\sum_{j\ge0}(1+\lambda_j)^s
       \langle f,e_j\rangle_{L^2(\mathrm{vol}_g)}^2,
 \qquad
 \|\eta\|_{-s}^2=\sum_{j\ge0}(1+\lambda_j)^{-s}\eta(e_j)^2.
\]
Thus
$\langle\eta,\zeta\rangle_{-s}
=\sum_{j\ge0}(1+\lambda_j)^{-s}\eta(e_j)\zeta(e_j)$.
For $\eta\in H^{-s}$ and $f\in H^s$, their dual pairing is $\eta(f)$.
The symbol $A^{\dagger_{-s}}$ denotes the Hilbert adjoint of a bounded
operator $A$ on $H^{-s}$, characterized by
$\langle A\eta,\zeta\rangle_{-s}
=\langle\eta,A^{\dagger_{-s}}\zeta\rangle_{-s}$.
For a nonnegative integer $r$, set
$\|f\|_{C^r}=\max_{0\le\ell\le r}\|\nabla^\ell f\|_\infty$,
where $\nabla^\ell$ is the $\ell$th covariant derivative,
$\nabla^0f=f$, and tensor norms use $g$ or the product metric on $M\times M$.

For a bounded nonnegative measurable $a:M\to\R$, write
\begin{equation}
 \mathcal V_\mu(a)=\int_M a(x)\,\delta_x\otimes\delta_x\,\mu(\dd x)
 \quad\hbox{on }H^{-s}.
 \label{eq:measure-covariance-notation}
\end{equation}
The integral is taken in trace norm.  This operator gives the covariance form $\int_M afg\,\dd\mu$ on test functions $f,g\in H^s$;
its trace-class construction is justified in Subsection~\ref{sec:compact-analysis}.

\begin{assumption}[Smooth kernels and type convergence]\label{ass:C}
The types and kernels are deterministic.  The kernels
$\kappa_n\in C^{2s}(M\times M)$ are symmetric and nonnegative.
There are symmetric kernels $\kappa,\Lambda\in C^{2s}(M\times M)$
and $\beta\in H^{-s}$ such that
\begin{gather}
 \|T_\kappa\|_{L^2(\mu)\to L^2(\mu)}>1,
 \label{eq:compact-model}\\
 \sqrt n(\mu^{(n)}-\mu)\longrightarrow\beta\quad\hbox{in }H^{-s},
 \qquad
 \sqrt n(\kappa_n-\kappa)\longrightarrow\Lambda\quad\hbox{in }C^{2s}.
 \label{eq:compact-expansions}
\end{gather}
The limiting kernel is irreducible with respect to $\mu$: if a
measurable $B\subset M$ has $\kappa=0$ $\mu\boxtimes\mu$-almost everywhere
on $B\times B^c$, then $\mu(B)\in\{0,1\}$.
\end{assumption}

A particle of type $x$ in the associated Poisson branching process has
offspring type intensity $\kappa(x,y)\mu(\dd y)$.  Its survival probability
$\rho$ is the maximal $[0,1]$-valued solution of
\begin{equation}
 \rho(x)=1-\exp\{-(\mathsf K(\rho\mu))(x)\}.
 \label{eq:compact-fixed-point}
\end{equation}
We use its smooth representative, satisfying this equation everywhere;
its construction is given in Subsection~\ref{sec:compact-analysis}.
Bollob\'as, Janson, and Riordan
\cite[Theorem~9.10]{BollobasJansonRiordan2007} give, under \assref{C},
\[
 \frac{\mathcal L_n(f)}n\xrightarrow{\mathrm p}\int_M f\rho\,\dd\mu,
 \qquad f\in C(M).
\]
In particular, $|C_1^{(n)}|/n\xrightarrow{\mathrm p}\int_M\rho\,\dd\mu>0$.
The derivative governing the fluctuations is
\begin{equation}
 J\eta=\eta-(1-\rho)(\mathsf K\eta)\mu,\qquad\eta\in H^{-s}.
 \label{eq:compact-J}
\end{equation}

\begin{theorem}[Gaussian limit for the type measure]\label{thm:compact}
Assume \assref{C}.  Then $J$ is boundedly invertible on $H^{-s}$ and,
in that space,
\begin{equation}
 \begin{aligned}
 \frac{\mathcal L_n-n\rho\mu}{\sqrt n}
 &\Rightarrow\Normal\!\Bigl(
 J^{-1}\{\rho\beta+(1-\rho)(\mathsf K_\Lambda(\rho\mu))\mu\},\\
 &\qquad J^{-1}\mathcal V_\mu(\rho(1-\rho))
                    (J^{-1})^{\dagger_{-s}}\Bigr).
 \end{aligned}
 \label{eq:compact-CLT}
\end{equation}
\end{theorem}

Here $\rho\beta+(1-\rho)(\mathsf K_\Lambda(\rho\mu))\mu$ is the shift
caused by the type and kernel perturbations, and
$\mathcal V_\mu(\rho(1-\rho))$ is the covariance of the edge noise.
As in Theorem~\ref{thm:countable-concrete}, $J$ is the derivative at
$\rho\mu$ of the mass balance map $h$ defined in
\eqref{eq:exploration-residuals}.  Its inverse transforms both
contributions into fluctuations of the giant type measure.

\paragraph{Remark (Scalar observations).}\phantomsection\label{rem:compact-scalar-observations}
Let $J^*:H^s\to H^s$ be the transpose under the dual pairing,
characterized by $(J\eta)(f)=\eta(J^*f)$.  It is given by
\begin{equation}
 (J^*f)(y)
 =f(y)-\int_M(1-\rho(x))\kappa(x,y)f(x)\,\mu(\dd x).
 \label{eq:compact-transpose}
\end{equation}
Theorem~\ref{thm:compact} gives bounded invertibility of $J$, hence of
$J^*$; write $R=(J^*)^{-1}$.  For $f,g\in H^s$, the limiting mean at
$f$ and the covariance at $f,g$ in that theorem are, respectively,
\begin{gather}
 \beta(\rho Rf)
 +\int_M (Rf)(1-\rho)(\mathsf K_\Lambda(\rho\mu))\,\dd\mu,
 \label{eq:compact-scalar-mean}\\
 \int_M (Rf)(Rg)\rho(1-\rho)\,\dd\mu.
 \label{eq:compact-scalar-covariance}
\end{gather}

Independent type sampling adds the empirical type fluctuation to the
edge noise, giving the following covariance.

\begin{corollary}[Independent type sampling]\label{cor:compact-iid}
Suppose that $\kappa\in C^{2s}(M\times M)$ is symmetric, nonnegative,
irreducible with respect to $\mu$, and supercritical as in
\eqref{eq:compact-model}.  Let $X_1,X_2,\ldots$ be independent with law
$\mu$, and construct the graph with types $X_1,\ldots,X_n$ and the fixed
kernel $\kappa$.  Then, in $H^{-s}$,
\begin{equation}
 \frac{\mathcal L_n-n\rho\mu}{\sqrt n}
 \Rightarrow\Normal\!\left(0,
 J^{-1}\{\mathcal V_\mu(\rho)-(\rho\mu)\otimes(\rho\mu)\}
 (J^{-1})^{\dagger_{-s}}\right).
 \label{eq:compact-iid-limit}
\end{equation}
\end{corollary}

\paragraph{Remark (Joint observables and type proportions).}\phantomsection\label{rem:type-proportions}
Suppose a main theorem or a corollary for independently sampled types
gives the limit $\Normal(m,C)$ on the Hilbert space $\mathcal H$ of the
model.  For countable types, take $f=(f_i)_i$ with
$\sum_i f_i^2/w_i<\infty$; this includes every bounded $f$ by
\eqref{eq:main-countable-summability}.  On a manifold, take $f\in H^s$.
These classes define continuous linear observations in the two spaces.
For a finite collection $f_1,\ldots,f_k$ in these classes, let
$T:\mathcal H\to\R^k$ collect the corresponding observations.
The continuous mapping theorem gives
\[
 \left(\frac{\mathcal L_n(f_j)-n(\rho\mu)(f_j)}{\sqrt n}\right)_{j=1}^k
 \Rightarrow\Normal_k(Tm,TCT^\dagger).
\]

For a single observable $f$, put
$a_f=(\rho\mu)(f)/(\rho\mu)(1)$.  Let $u_f\in\mathcal H$ be the Riesz
representative of the continuous functional
\[
 \eta\longmapsto\frac{\eta(f)-a_f\eta(1)}{(\rho\mu)(1)},
\]
using the mass-vector identification for countable types.
Since $(\rho\mu)(1)>0$, the delta method yields
\begin{equation}
 \sqrt n\left\{\frac{\mathcal L_n(f)}{\mathcal L_n(1)}-a_f\right\}
 \Rightarrow\Normal_1\!\left(
 \langle m,u_f\rangle_{\mathcal H},
 \langle Cu_f,u_f\rangle_{\mathcal H}\right).
 \label{eq:giant-type-proportions}
\end{equation}
The ratio may be defined arbitrarily on $\{\mathcal L_n(1)=0\}$,
whose probability tends to zero.  Thus the same Gaussian limits describe
averages of type observables within the giant, jointly with its size.

\begin{example}[Heat kernels and spectral fluctuations]\label{ex:compact-smooth}
Let $d\ge1$, $\mu=\mathrm{vol}_g$, $u>0$, $\lambda>1$, and
$0<\varepsilon<\lambda$.  The heat kernel at time $u$ has the expansion
\begin{equation}
 p_u(x,y)=\sum_{j\ge0}e^{-u\lambda_j}e_j(x)e_j(y),\qquad
 \int_Mp_u(x,y)\,\mu(\dd y)=1.
 \label{eq:heat-kernel-expansion}
\end{equation}
Take the fixed kernel
\begin{equation}
 \kappa(x,y)=\lambda+\varepsilon\{p_u(x,y)-1\}.
 \label{eq:heat-kernel-model}
\end{equation}
It is smooth, strictly positive, supercritical, and of infinite rank.
Its row integral is $\lambda$, so $\rho$ is the positive constant solving
$\rho=1-e^{-\lambda\rho}$.  Put
\begin{equation}
 \theta_0=\lambda,\qquad
 \theta_j=\varepsilon e^{-u\lambda_j}\ (j\ge1),\qquad
 r_j=1-(1-\rho)\theta_j\ (j\ge0).
 \label{eq:heat-kernel-modes}
\end{equation}
Then $T_\kappa e_j=\theta_je_j$ and
$r_j\ge1-(1-\rho)\lambda>0$.
For iid types, define
\[
 \mathcal C_{\mathrm{iid}}
 =\sum_{j\ge0}\sigma_j^2(e_j\mu)\otimes(e_j\mu),\qquad
 \sigma_j^2=
 \begin{cases}
 \rho(1-\rho)/r_0^2,&j=0,\\
 \rho/r_j^2,&j\ge1.
 \end{cases}
\]
Corollary~\ref{cor:compact-iid} gives
\begin{equation}
 \frac{\mathcal L_n-n\rho\mu}{\sqrt n}
 \Rightarrow\Normal(0,\mathcal C_{\mathrm{iid}})
 \quad\hbox{in }H^{-s}.
 \label{eq:heat-iid-modes}
\end{equation}
Equivalently, for every $m\ge0$,
\[
 \left(\frac{\mathcal L_n(e_j)-n\rho\ind_{\{j=0\}}}{\sqrt n}\right)_{j=0}^m
 \Rightarrow\Normal_{m+1}\!\left(0,\diag(\sigma_0^2,\ldots,\sigma_m^2)\right).
\]
Thus the limiting spectral coordinates are independent, and coordinate
$j$ has variance $\sigma_j^2$.

For deterministic types with
$\sqrt n\,\|\mu^{(n)}-\mu\|_{-s}\to0$, define
\[
 \mathcal C_{\mathrm{det}}
 =\sum_{j\ge0}\frac{\rho(1-\rho)}{r_j^2}(e_j\mu)\otimes(e_j\mu).
\]
Theorem~\ref{thm:compact} gives
\begin{equation}
 \frac{\mathcal L_n-n\rho\mu}{\sqrt n}
 \Rightarrow\Normal(0,\mathcal C_{\mathrm{det}})
 \quad\hbox{in }H^{-s}.
 \label{eq:heat-deterministic-modes}
\end{equation}
Both covariance series converge in trace norm, as verified in
Subsection~\ref{sec:compact-heat-proof}.  On a flat torus, regular grids
with $n=m^d$ points satisfy the deterministic type condition.
The spectral coordinate variances are $\sigma_j^2$ for iid types and
$\rho(1-\rho)/r_j^2$ for deterministic types.  In both cases,
\[
 \frac{|C_1^{(n)}|-n\rho}{\sqrt n}
 \Rightarrow\Normal_1\!\left(0,
 \frac{\rho(1-\rho)}{(1-(1-\rho)\lambda)^2}\right).
\]
Thus total size has the same limiting variance as for the constant kernel
$\lambda$, while the other coordinates depend on the Laplace spectrum.

For the deterministic types above, take
\begin{equation}
 \Lambda(x,y)=e_1(x)+e_1(y),\qquad
 \kappa_n=\kappa+n^{-1/2}\Lambda
 \label{eq:heat-kernel-perturbation}
\end{equation}
for all sufficiently large $n$.  The perturbation produces a nonzero
limiting mean:
\begin{equation}
 \frac{\mathcal L_n-n\rho\mu}{\sqrt n}
 \Rightarrow\Normal\!\left(
 \frac{\rho(1-\rho)}{r_1}e_1\mu,\mathcal C_{\mathrm{det}}\right)
 \quad\hbox{in }H^{-s}.
 \label{eq:heat-kernel-perturbed-limit}
\end{equation}
The total-size limiting mean remains zero, while the mean of the first
spectral coordinate is $\rho(1-\rho)/r_1$.
\end{example}

\section{Exploration and local linearization}\label{sec:proof}

A vertex exploration started from a small random seed set approximates
the giant component.  We describe this algorithm first, then compensate
the loss of unexplored vertices to obtain a martingale.  Conservation of
the type measure gives an exact equation relating the explored mass to
its martingale error.  Finally, a local linearization lemma converts this
equation into a first order fluctuation formula: the limiting derivative
transforms the exploration noise and the parameter perturbations into
the fluctuation of the giant.

\subsection{The seeded exploration}\label{sec:exploration}

We explore all components meeting a small random seed set; this reaches
the giant with high probability while adding a negligible number of other
vertices.  The unexplored vertices will provide a martingale description
of the edge randomness.

We use a version of the classical vertex exploration, started from a
random set of seeds.  For Erd\H{o}s--R\'enyi graphs, the exploration and
its martingale analysis are described in
\cite[Section~2]{BollobasRiordan2012}; related epidemic sampling
procedures appear in Martin-L\"of \cite{MartinLof1986}.

We use the environment $\mathcal E_n$ from Section~\ref{sec:model}.
For the estimates in this section suppose that
$0\leq\kappa_n\leq\kappa_*$ with the deterministic bound specified there.

Put $\ell_n=\lfloor\log n\rfloor$.  We need the number of seeds to tend
to infinity in order to reach a component of positive limiting density.
The component estimates proved below give $|C_2^{(n)}|=\Op(\log n)$,
so $\ell_n\log n=o(\sqrt n)$ makes the number of extra explored vertices
negligible on the fluctuation scale.  Our choice satisfies both conditions;
the weighted error in the countable model is checked separately in
Lemma~\ref{lem:countable-seed-approximation}.
Independently of the environment and all edge variables, choose a
uniform subset $I_{n,0}\subset[n]$ of size $\ell_n$.  At time $t$, partition
the vertices into removed, active, and unexplored sets $R_{n,t},I_{n,t},U_{n,t}$.
Initially $R_{n,0}=\varnothing$, the active set is the seed set $I_{n,0}$, and
$U_{n,0}=[n]\setminus I_{n,0}$.  If
$I_{n,t}\ne\varnothing$, choose an active vertex $\omega_{n,t}$ by a fixed
measurable rule depending only on the exploration history, reveal its edges
to $U_{n,t}$, and put
\[
 B_{n,t}=\{v\in U_{n,t}:\omega_{n,t}\sim_n v\}.
\]
Set
\[
 R_{n,t+1}=R_{n,t}\cup\{\omega_{n,t}\},\quad
 I_{n,t+1}=(I_{n,t}\setminus\{\omega_{n,t}\})\cup B_{n,t},\quad
 U_{n,t+1}=U_{n,t}\setminus B_{n,t}.
\]
Let $\tau_n=\inf\{t\geq0:I_{n,t}=\varnothing\}$ and extend all sets
constantly after $\tau_n$.  Each step removes one vertex, so $\tau_n\leq n$.
Let $\mathcal F_{n,t}$ be the sigma algebra generated by the environment,
the seeds, and all choices and edge indicators revealed before step $t$.  Then $\tau_n$ is a stopping time
and the choice of $\omega_{n,t}$ is $\mathcal F_{n,t}$ measurable.
At termination, $R_{n,\tau_n}$ is the union of the components meeting $I_{n,0}$.

Two estimates explain the choice of seeds.  Conditional on the graph,
the seed set is uniform among the $\ell_n$ element subsets of $[n]$.
Sampling without replacement therefore gives
\begin{equation}
 \Pp(I_{n,0}\cap C_1^{(n)}=\varnothing\mid G_n)
 \leq\left(1-\frac{|C_1^{(n)}|}{n}\right)^{\ell_n},
 \qquad
 |R_{n,\tau_n}\setminus C_1^{(n)}|\leq\ell_n|C_2^{(n)}|.
 \label{eq:seed-approximation}
\end{equation}
For the first inequality, the exact probability is
$\binom{n-|C_1^{(n)}|}{\ell_n}/\binom n{\ell_n}$, with value zero if
$n-|C_1^{(n)}|<\ell_n$.  In particular, for every fixed $a\in(0,1)$,
\begin{equation}
 \Pp(I_{n,0}\cap C_1^{(n)}=\varnothing)
 \leq\Pp(|C_1^{(n)}|/n<a)+(1-a)^{\ell_n}.
 \label{eq:seed-hit}
\end{equation}
The two terms tend to zero whenever the giant has a limiting density
larger than $a$.  The second inequality in \eqref{eq:seed-approximation}
holds on the event that the seeds meet the giant: each additional
component contains a seed.  Together with $|C_2^{(n)}|=\Op(\log n)$,
it bounds the number of extra explored vertices by
$\Op((\log n)^2)$.
Lemmas~\ref{lem:countable-seed-approximation}
and~\ref{lem:compact-seeded-approximation} bound the type counts or
measures of these extra vertices in the respective norms.

We keep track of types rather than only the number of explored vertices.
Define the finite atomic measures
\begin{equation}
 \mathcal R_{n,t}=\sum_{v\in R_{n,t}}\delta_{x_v^{(n)}},\qquad
 \mathcal U_{n,t}=\sum_{v\in U_{n,t}}\delta_{x_v^{(n)}},\qquad
 S_n=\sum_{v\in I_{n,0}}\delta_{x_v^{(n)}}.
 \label{eq:exploration-measures}
\end{equation}
At termination there are no active vertices, so
\begin{equation}
 \mathcal R_{n,\tau_n}+\mathcal U_{n,\tau_n}=n\mu^{(n)}.
 \label{eq:measure-conservation}
\end{equation}
Thus information about the remaining unexplored measure also determines
the explored measure, which approximates $\mathcal L_n$.

The measure $\mathcal U_{n,t}$ decreases at each step: an unexplored vertex
of type $x$ remains unexplored with conditional probability
$1-p_n(x,x_{\omega_{n,t}}^{(n)})$.  We denote the product of these one step
survival factors by
\begin{equation}
 Q_{n,t}(x)=\prod_{s<t\wedge\tau_n}
            \{1-p_n(x,x_{\omega_{n,s}}^{(n)})\}
 =\exp\left\{-\frac1n\sum_{v\in R_{n,t\wedge\tau_n}}
                          \kappa_n(x,x_v^{(n)})\right\}.
 \label{eq:exploration-Q}
\end{equation}
Dividing the remaining mass by this product compensates its predictable
loss.  This motivates
\begin{equation}
 W_{n,t}=Q_{n,t}^{-1}\mathcal U_{n,t},\qquad
 M_{n,t}=W_{n,t}/n.
 \label{eq:exploration-W-M}
\end{equation}
With the measure notation of Subsection~\ref{sec:notation},
\eqref{eq:exploration-W-M} gives
\[
 W_{n,t}(f)=\sum_{v\in U_{n,t}}\frac{f(x_v^{(n)})}{Q_{n,t}(x_v^{(n)})}.
\]
In particular, $W_{n,0}=n\mu^{(n)}-S_n$ and
$e^{-\kappa_*}\leq Q_{n,t}\leq1$.  Lemma~\ref{lem:exploration-martingales}
shows that $W_{n,\cdot}(f)$ is a martingale.  Its fluctuation at $\tau_n$
gives the edge noise, while $M_{n,t}$ stays close to $\mu^{(n)}$ and allows
us to identify the limiting covariance.

We next express the exploration in a form suitable for a fluctuation
limit.  Normalize the explored measure by the graph size and the
martingale increment by the fluctuation scale:

\begin{equation}
 \widehat\mu_n=\mathcal R_{n,\tau_n}/n,\qquad
 N_n=\frac{W_{n,\tau_n}-W_{n,0}}{\sqrt n}.
 \label{eq:exploration-martingale-noise}
\end{equation}
The unexplored measure differs from $nQ_{n,\tau_n}\mu^{(n)}$ by
\begin{equation}
 Z_n=\frac{\mathcal U_{n,\tau_n}-nQ_{n,\tau_n}\mu^{(n)}}{\sqrt n}
     =Q_{n,\tau_n}N_n-Q_{n,\tau_n}S_n/\sqrt n.
 \label{eq:exploration-noise}
\end{equation}
By \eqref{eq:exploration-Q}, the survival factor at the stopping time is
$e^{-\mathsf K_{\kappa_n}\widehat\mu_n}$.  Thus conservation
\eqref{eq:measure-conservation} compares the explored mass with the
remaining mass predicted by this factor.  The discrepancy is described
by the following maps on finite signed measures:

\begin{equation}
 h_n(\eta)=\eta+e^{-\mathsf K_{\kappa_n}\eta}\mu^{(n)}-\mu^{(n)},
 \qquad h(\eta)=\eta+e^{-\mathsf K_\kappa\eta}\mu-\mu.
 \label{eq:exploration-residuals}
\end{equation}
The limiting fixed point equation gives $h(\rho\mu)=0$.  We will
therefore expand the exploration equation near $\rho\mu$; its derivative
will determine how errors in mass balance affect the explored mass.
The bounded kernels make the displayed expressions finite signed
measures.  Lemmas~\ref{lem:countable-terminal-maps}
and~\ref{lem:compact-terminal-maps} extend these maps to the respective
Hilbert spaces and establish their differentiability.

\begin{lemma}[The exploration equation]\label{lem:exploration-equation}
The normalized explored measure satisfies the exact identity
\begin{equation}
 \sqrt n\,h_n(\widehat\mu_n)+Z_n=0.
 \label{eq:exploration-equation}
\end{equation}
\end{lemma}

\begin{proof}
By \eqref{eq:exploration-Q},
$Q_{n,\tau_n}=e^{-\mathsf K_{\kappa_n}\widehat\mu_n}$.
Substitute this into \eqref{eq:exploration-residuals} and use conservation
\eqref{eq:measure-conservation}:
\[
 \sqrt n\,h_n(\widehat\mu_n)
 =\frac{\mathcal R_{n,\tau_n}+nQ_{n,\tau_n}\mu^{(n)}-n\mu^{(n)}}{\sqrt n}
 =-Z_n.
\]
\end{proof}

\subsection{Martingale identities and finite projections}
\label{sec:stopped-martingales}

The martingale central limit theorem requires a limiting conditional
covariance and control of large increments. The increment identities
below provide the covariance; their fourth moment bounds give the
Lindeberg condition in both models.

For $t<\tau_n$ and $v\in U_{n,t}$, let $B_{n,t,v}$ be the indicator of the edge
from $\omega_{n,t}$ to $v$, and put
$p_{n,t}(x)=p_n(x,x_{\omega_{n,t}}^{(n)})$.
Write $W_{n,\cdot}(f)=(W_{n,t}(f))_{0\leq t\leq n}$ for the whole
process.  For two real square integrable $(\mathcal F_{n,t})$ martingales
$X_{n,\cdot}$ and $Y_{n,\cdot}$, their predictable covariation, the
accumulated conditional covariance of their increments, is
\[
 \langle X_{n,\cdot},Y_{n,\cdot}\rangle_t
 =\sum_{s<t}\E(\Delta X_{n,s}\Delta Y_{n,s}\mid\mathcal F_{n,s}).
\]

We write $\langle X_{n,\cdot}\rangle_t
=\langle X_{n,\cdot},X_{n,\cdot}\rangle_t$ for predictable quadratic
variation.

\begin{lemma}[Exploration martingales]\label{lem:exploration-martingales}
For every bounded measurable $f:\mathcal S\to\R$, the process
$(W_{n,t}(f):0\leq t\leq n)$ is a square integrable
$(\mathcal F_{n,t})$ martingale, constant after $\tau_n$.
For $t<\tau_n$, its increment is
\begin{equation}
 \Delta W_{n,t}(f)=\sum_{v\in U_{n,t}}
 \frac{f(x_v^{(n)})\{p_{n,t}(x_v^{(n)})-B_{n,t,v}\}}
      {Q_{n,t}(x_v^{(n)})(1-p_{n,t}(x_v^{(n)}))}.
 \label{eq:measure-increment}
\end{equation}
For two such functions $f,g$,
\begin{equation}
 \frac1n\langle W_{n,\cdot}(f),W_{n,\cdot}(g)\rangle_{\tau_n}
 =\sum_{t<\tau_n}M_{n,t}\bigl(fg(Q_{n,t+1}^{-1}-Q_{n,t}^{-1})\bigr).
 \label{eq:common-bracket}
\end{equation}
\end{lemma}

\begin{proof}
The measures $\mathcal U_{n,t}$ and functions $Q_{n,t}$ are
$\mathcal F_{n,t}$ measurable, so $W_{n,t}(f)$ is adapted.  Since
$e^{-\kappa_*}\leq Q_{n,t}\leq1$ and $|U_{n,t}|\leq n$,
\[
 |W_{n,t}(f)|
 =\left|\sum_{v\in U_{n,t}}\frac{f(x_v^{(n)})}{Q_{n,t}(x_v^{(n)})}\right|
 \leq ne^{\kappa_*}\|f\|_\infty.
\]
Thus $W_{n,t}(f)$ belongs to $L^2$ for every $t$ and $n$.

Conditional on $\mathcal E_n$, the edge indicators are independent,
and the seeds are independent of them.  The exploration history depends
only on the seeds and previously revealed edges.  Since the choice of
$\omega_{n,t}$ is $\mathcal F_{n,t}$ measurable, unrevealed edges retain their
conditionally independent Bernoulli laws given $\mathcal F_{n,t}$.
For $t<\tau_n$, the edges $\{\omega_{n,t},v\}$, $v\in U_{n,t}$, are distinct and
unrevealed: every earlier test involved a vertex in $R_{n,t}$, whereas neither
endpoint of a current test lies in $R_{n,t}$.  Consequently, given
$\mathcal F_{n,t}$, the variables $B_{n,t,v}$ are independent with means
$p_{n,t}(x_v^{(n)})$.

The exploration updates give
\[
 Q_{n,t+1}(x)=Q_{n,t}(x)(1-p_{n,t}(x)),\qquad
 \mathcal U_{n,t+1}=\sum_{v\in U_{n,t}}(1-B_{n,t,v})\delta_{x_v^{(n)}}.
\]
It follows that
\[
 W_{n,t+1}(f)-W_{n,t}(f)
 =\sum_{v\in U_{n,t}}\frac{f(x_v^{(n)})}{Q_{n,t}(x_v^{(n)})}
   \left\{\frac{1-B_{n,t,v}}{1-p_{n,t}(x_v^{(n)})}-1\right\},
\]
which is \eqref{eq:measure-increment}.  All its coefficients are
$\mathcal F_{n,t}$ measurable, and
$\E[B_{n,t,v}\mid\mathcal F_{n,t}]=p_{n,t}(x_v^{(n)})$.  Hence
$\E[\Delta W_{n,t}(f)\mid\mathcal F_{n,t}]=0$ for $t<\tau_n$.
For $t\geq\tau_n$, the increment is zero by definition.  Together with
adaptedness and square integrability, this proves the martingale assertion.

Conditional independence also gives
\[
 \E[\Delta W_{n,t}(f)\Delta W_{n,t}(g)\mid\mathcal F_{n,t}]
 =\sum_{v\in U_{n,t}}
 \frac{f(x_v^{(n)})g(x_v^{(n)})p_{n,t}(x_v^{(n)})}
      {Q_{n,t}(x_v^{(n)})^2(1-p_{n,t}(x_v^{(n)}))}.
\]
Since
\[
 Q_{n,t+1}^{-1}-Q_{n,t}^{-1}=\frac{p_{n,t}}{Q_{n,t}(1-p_{n,t})},\qquad
 M_{n,t}=\frac{\mathcal U_{n,t}}{nQ_{n,t}},
\]
division by $n$ and summation over $t<\tau_n$ prove
\eqref{eq:common-bracket}.

\end{proof}

The covariance in \eqref{eq:common-bracket} determines the Gaussian
limit once $M_{n,t}$ and $Q_{n,\tau_n}$ are understood.  To apply a
martingale central limit theorem, we also need large increments to be
negligible.  Lemma~\ref{lem:exploration-lindeberg} gives this conditional
Lindeberg condition in \eqref{eq:common-lindeberg}.

\begin{lemma}[Bounds for large martingale increments]
\label{lem:exploration-lindeberg}
For every bounded measurable $f:\mathcal S\to\R$, every
$0\leq t<n$, and every $\varepsilon>0$,
\begin{equation}
 \begin{gathered}
 \E\bigl((\Delta W_{n,t}(f))^4\mid\mathcal F_{n,t}\bigr)
 \leq C\|f\|_\infty^4,\\
 \sum_{t<n}\E\left[
 \frac{(\Delta W_{n,t}(f))^2}{n}
 \ind_{\{|\Delta W_{n,t}(f)|>\varepsilon\sqrt n\}}
 \biggm|\mathcal F_{n,t}\right]
 \leq\frac{C\|f\|_\infty^4}{\varepsilon^2n},
 \end{gathered}
 \label{eq:common-lindeberg}
\end{equation}
where $C$ depends only on $\kappa_*$.
\end{lemma}

\begin{proof}
To bound the fourth moments, let $B_v$ be independent Bernoulli variables
of means $p_v$, put $X_v=B_v-p_v$, and let $b_v$ be real coefficients.
Centering and independence give
\begin{align*}
 \E\left(\sum_v b_vX_v\right)^4
 &=\sum_v b_v^4\E X_v^4
   +6\sum_{v<u}b_v^2b_u^2\E X_v^2\E X_u^2\\
 &\leq\sum_v b_v^4p_v
   +3\left(\sum_v b_v^2p_v\right)^2,
\end{align*}
where $\E X_v^2\leq p_v$ and $\E X_v^4\leq p_v$.
Apply this calculation conditionally on $\mathcal F_{n,t}$ to
\eqref{eq:measure-increment}.  Its coefficients satisfy
\[
 |b_v|=\frac{|f(x_v^{(n)})|}{Q_{n,t+1}(x_v^{(n)})}
 \leq e^{\kappa_*}\|f\|_\infty,\qquad
 \sum_{v\in U_{n,t}}p_{n,t}(x_v^{(n)})\leq n\frac{\kappa_*}{n}=\kappa_*.
\]
This proves the fourth moment bound.  Finally,
\[
 \frac{(\Delta W_{n,t}(f))^2}{n}
 \ind_{\{|\Delta W_{n,t}(f)|>\varepsilon\sqrt n\}}
 \leq\frac{(\Delta W_{n,t}(f))^4}{\varepsilon^2n^2}.
\]
Taking conditional expectations and summing over at most $n$ steps proves
the second bound in \eqref{eq:common-lindeberg}.
\end{proof}

Finite projections alone do not give convergence in a Hilbert space:
the omitted coordinates must also be controlled.  We use the following
criterion for coordinate projections, Laplace projections, and their
products.

\begin{lemma}[Convergence from finite projections]\label{lem:hilbert-projections}
Let $\mathcal H$ be a separable real Hilbert space and let $P_m$ be
finite rank orthogonal projections such that $P_mx\to x$ for every
$x\in\mathcal H$.  Let $Y_n,Y$ be $\mathcal H$ valued random elements.
Suppose that $P_mY_n\Rightarrow P_mY$ for every fixed $m$, and that,
for every $\varepsilon>0$,
\begin{equation}
 \lim_{m\to\infty}\limsup_{n\to\infty}
 \Pp(\|(I-P_m)Y_n\|>\varepsilon)=0.
 \label{eq:hilbert-projection-tails}
\end{equation}
Then $Y_n\Rightarrow Y$ in $\mathcal H$.
Condition~\eqref{eq:hilbert-projection-tails} follows in particular from
$\lim_{m\to\infty}\sup_n\E\|(I-P_m)Y_n\|^2=0$, or from tightness of
$(Y_n)_n$.
\end{lemma}

\begin{proof}
For every function $F$ with $\|F\|_\infty\leq1$ and Lipschitz constant
at most one,
\[
 |\E F(Y_n)-\E F(P_mY_n)|
 \leq\varepsilon+2\Pp(\|(I-P_m)Y_n\|>\varepsilon).
\]
The analogous difference for $Y$ tends to zero as $m\to\infty$, since
$P_mY\to Y$ almost surely.  For fixed $m$, finite dimensional convergence
gives $\E F(P_mY_n)\to\E F(P_mY)$.  First let $n\to\infty$, then
$m\to\infty$, and finally $\varepsilon\downarrow0$; convergence against
bounded Lipschitz functions proves the claim.
The second moment condition implies \eqref{eq:hilbert-projection-tails}
by Markov's inequality.  If $(Y_n)_n$ is tight, choose a compact set
carrying arbitrarily large probability uniformly in $n$.  Strong
convergence of the orthogonal projections is uniform on this compact
set, because $\|I-P_m\|\leq1$; hence
\eqref{eq:hilbert-projection-tails} holds again.
\end{proof}

\subsection{Local linearization}\label{sec:terminal-linearization}

The exploration equation \eqref{eq:exploration-equation} gives the
random error in mass balance.  To recover the fluctuation of the explored
mass, we expand its left side at the limiting giant and invert the
derivative.  The next lemma justifies this step and proves the
$n^{-1/2}$ rate required to make the Taylor remainder negligible.

Let $\mathcal H$ be a real Banach space with norm $\|\cdot\|$, fix a
deterministic point $x_\star\in\mathcal H$, and let $\mathcal U$ be an open
neighborhood of $x_\star$.  Let $g:\mathcal U\to\mathcal H$ be deterministic
and continuously Fr\'echet differentiable, and write $J=Dg(x_\star)$.
On a probability space $(\Omega,\mathcal F,\Pp)$, let
$g_n:\Omega\times\mathcal U\to\mathcal H$ be jointly measurable maps
such that $g_n(\omega,\cdot)$ is continuously Fr\'echet differentiable
almost surely.  We suppress $\omega$ in the notation; derivatives act on the spatial
variable, and operator norms are understood as in
Subsection~\ref{sec:notation}.
Let $X_n$ be $\mathcal H$ valued random
elements on the same probability space; set $g_n(X_n)=0$ when
$X_n\notin\mathcal U$.

In our applications, sampled types or random kernels make $g_n$ random,
so $g_n(x_\star)$ and the derivative supremum below may be random even though
$x_\star$ is fixed.  When $g_n$ is deterministic, the $\Op(1)$ bound on
$\sqrt n\,g_n(x_\star)$ reduces to an $O(1)$ bound, and the derivative
convergence is ordinary convergence.  The value $g_n(X_n)$ may still
be random through $X_n$.

\begin{lemma}[Quantitative local linearization]\label{lem:linearization}
Suppose that $g(x_\star)=0$, $J$ is boundedly invertible, and
$X_n\xrightarrow{\mathrm p}x_\star$.
Assume
\[
 \sqrt n\,g_n(x_\star)=\Op(1),\qquad
 \sqrt n\,g_n(X_n)=\Op(1),
\]
and suppose that there is $r_0>0$ such that the closed ball of radius $r_0$
about $x_\star$ lies in $\mathcal U$ and
\[
 \sup_{\|x-x_\star\|\leq r_0}\|Dg_n(x)-Dg(x)\|
 \xrightarrow{\mathrm p}0.
\]
Then
\[
 \sqrt n(X_n-x_\star)=\Op(1)
\]
and
\begin{equation}
 \sqrt n\,g_n(X_n)
 =\sqrt n\,g_n(x_\star)+J\sqrt n(X_n-x_\star)+\op(1).
 \label{eq:quantitative-linearization}
\end{equation}
\end{lemma}

\begin{proof}[Proof of Lemma~\ref{lem:linearization}]
We first obtain the $n^{-1/2}$ rate using uniform invertibility near $x_\star$;
only then do we replace the averaged derivative by $J$.
By continuity of $Dg$ at $x_\star$, choose $0<\delta\leq r_0$ so that
\[
 \sup_{\|x-x_\star\|\leq\delta}\|Dg(x)-J\|
 \leq \frac1{4\|J^{-1}\|}.
\]
For $\|x-x_\star\|\leq\delta$, define the Bochner integral of bounded operators
\[
 B_n(x)=\int_0^1Dg_n(x_\star+s(x-x_\star))\,\dd s.
\]
The fundamental theorem of calculus gives
$g_n(x)-g_n(x_\star)=B_n(x)(x-x_\star)$.
The assumed derivative convergence implies that, with probability tending
to one,
\[
 \sup_{\|x-x_\star\|\leq\delta}\|B_n(x)-J\|
 \leq \frac1{2\|J^{-1}\|}.
\]

On this event, \(\|J^{-1}(B_n(x)-J)\|\leq\frac12\) throughout the ball.
The Neumann series shows that every
$B_n(x)$ in the ball is boundedly invertible and
$\sup_{\|x-x_\star\|\leq\delta}\|B_n(x)^{-1}\|\leq2\|J^{-1}\|$.
Since $X_n$ lies in the ball with probability tending to one,
\[
 \sqrt n(X_n-x_\star)
 =B_n(X_n)^{-1}\sqrt n\{g_n(X_n)-g_n(x_\star)\}=\Op(1).
\]
Finally, on the same event and with $X_n$ in the ball,
\begin{align*}
 \|B_n(X_n)-J\|
 &\leq\sup_{\|x-x_\star\|\leq r_0}\|Dg_n(x)-Dg(x)\|\\
 &\quad+\sup_{0\leq t\leq1}
          \|Dg(x_\star+t(X_n-x_\star))-Dg(x_\star)\|=\op(1).
\end{align*}
The first term tends to zero by assumption, and the second by continuity
of $Dg$ at $x_\star$ and $X_n\xrightarrow{\mathrm p}x_\star$.
Multiplying this difference by
$\sqrt n(X_n-x_\star)=\Op(1)$ proves
\eqref{eq:quantitative-linearization}.
\end{proof}

For the graph applications, take $x_\star=\rho\mu$ and let $g_n$
be the extension of $h_n$ from \eqref{eq:exploration-residuals}.
The derivative is $J=Dh(\rho\mu)$.  The parameter expansion and the
martingale limit will give, for deterministic perturbations,
\[
 -\sqrt n\,h_n(\rho\mu)\longrightarrow b,\qquad
 Z_n\Rightarrow\Normal(0,\Sigma),
\]
where $b$ is the limiting parameter drift and $\Sigma$ is the covariance
operator of the exploration noise on the relevant Hilbert space.  Equation~\eqref{eq:exploration-equation} then gives
$\sqrt n\,g_n(\widehat\mu_n)=-Z_n=\Op(1)$.  Once the explored measure
converges to $\rho\mu$, the derivatives converge locally, and $J$ has a
bounded inverse, \eqref{eq:quantitative-linearization} yields
\begin{equation}
 \sqrt n(\widehat\mu_n-\rho\mu)=J^{-1}(b-Z_n)+\op(1)
 \Rightarrow\Normal\!\left(J^{-1}b,
 J^{-1}\Sigma(J^{-1})^\dagger\right).
 \label{eq:linearized-exploration}
\end{equation}
The following table identifies where these conditions and the
approximation of the giant by the exploration are verified.

\begin{center}
\small
\begin{tabular}{@{}p{0.36\textwidth}p{0.27\textwidth}p{0.27\textwidth}@{}}
\hline
Property & Countable types & Manifold types\\
\hline
$\widehat\mu_n\to \rho\mu$ and exploration error $o_{\mathrm p}(\sqrt n)$
 & Lemma~\ref{lem:countable-seed-approximation}
 & Lemma~\ref{lem:compact-seeded-approximation}\\[3pt]
Martingale limit and Gaussian law of $Z_n$
 & Lemmas~\ref{lem:countable-martingale-limit}, \ref{lem:countable-terminal-equation}
 & Lemmas~\ref{lem:compact-martingale-limit}, \ref{lem:compact-terminal-noise}\\[3pt]
$-\sqrt n\,h_n(\rho\mu)\to b$
 & Lemma~\ref{lem:countable-terminal-maps}
 & Lemma~\ref{lem:compact-terminal-maps}\\[3pt]
Local convergence of the derivatives
 & Lemma~\ref{lem:countable-terminal-maps}
 & Lemma~\ref{lem:compact-terminal-maps}\\[3pt]
Bounded inverse of $J$
 & Lemma~\ref{lem:countable-inverse}
 & Lemma~\ref{lem:compact-inverse}\\
\hline
\end{tabular}
\end{center}

For independently sampled types, the first order parameter term is
random.  Lemmas~\ref{lem:countable-iid-joint}
and~\ref{lem:compact-joint-noise} identify its joint limit with $Z_n$;
the same linearization then gives the fluctuation limit.

\section{Proofs for countably many types}\label{sec:countable-proof}

Weighted tail bounds control the types omitted by finite coordinate
projections.  Together with a bounded inverse of the fixed point derivative,
they turn the exploration equation \eqref{eq:exploration-equation} into
the type vector limit of Theorem~\ref{thm:countable-concrete}.

\subsection{Weighted convergence and the fixed point equation}\label{sec:countable-estimates}

Weighted moment convergence controls empirical tails and upgrades the
coordinatewise law of large numbers for the giant to convergence in
$\mathcal H_w$. The resulting estimates also give the differentiability
and invertibility needed to linearize $h_n$ in
\eqref{eq:exploration-equation}.

\begin{lemma}[Weighted moment criteria]\label{lem:countable-weighted-moments}
There exist weights $w_i\geq1$ satisfying
\eqref{eq:main-countable-summability} if and only if
$\sum_i\sqrt{\mu_i}<\infty$.
For fixed weights with $\sum_iw_i\mu_i<\infty$, let
$\mu^{(n)}$ be empirical probability vectors such that
$\mu_i^{(n)}\to\mu_i$ for every $i$.  Then the uniform weighted tail
condition in \eqref{eq:main-countable-counts} is equivalent to
\begin{equation}
 \sum_iw_i\mu_i^{(n)}\longrightarrow\sum_iw_i\mu_i.
 \label{eq:countable-weighted-moments}
\end{equation}
Either condition implies
$\sum_iw_i|\mu_i^{(n)}-\mu_i|\to0$.
\end{lemma}

\begin{proof}
Cauchy--Schwarz gives
\[
 \sum_i\sqrt{\mu_i}
 \leq\left(\sum_iw_i\mu_i\right)^{1/2}
       \left(\sum_iw_i^{-1}\right)^{1/2}.
\]
Conversely, if $\sum_i\sqrt{\mu_i}<\infty$, then $w_i=\mu_i^{-1/2}\geq1$
satisfies both conditions in \eqref{eq:main-countable-summability}.

For the weighted moment criterion, assume
\eqref{eq:countable-weighted-moments}.  Dominated convergence gives
\[
 \sum_iw_i\min(\mu_i^{(n)},\mu_i)
 \longrightarrow\sum_iw_i\mu_i,
\]
because the summands are bounded by $w_i\mu_i$.
Thus
\[
 \sum_iw_i|\mu_i^{(n)}-\mu_i|\longrightarrow0.
\]
For all sufficiently large $n$, this bounds the empirical tail by the
limiting tail plus an arbitrarily small error.  The finitely many remaining
empirical vectors have finite support, so their tails also vanish.
This proves the uniform tail condition.  Conversely, split the weighted
sums at a fixed $m$, use coordinate convergence on $\mathcal I_m$, and
then let $m\to\infty$; the two tails vanish uniformly.

\end{proof}

For the remainder of this subsection, assume the countable type model
and \eqref{eq:main-countable-summability}, with deterministic types and
uniformly bounded kernels.

\begin{lemma}[Giant component estimates]
\label{lem:countable-components}
Suppose that
\[
 \mu^{(n)}\longrightarrow\mu\quad\hbox{in }\mathcal H_w,\qquad
 \lim_{m\to\infty}\sup_n\sum_{i>m}w_i\mu_i^{(n)}=0,\qquad
 \|\Delta_n\|_{w,\infty}\longrightarrow0.
\]
Then
\begin{equation}
 \left\|\frac{L_n}{n}-\mu\rho\right\|_1+
 \left\|\frac{L_n}{n}-\mu\rho\right\|_w\xrightarrow{\mathrm p}0,
 \qquad |C_2^{(n)}|=\Op(\log n).
 \label{eq:countable-first-order}
\end{equation}
In particular, $|C_1^{(n)}|/n\xrightarrow{\mathrm p}\bar\rho$.
\end{lemma}

\begin{proof}
This is a consequence of
\cite[Theorems~9.10 and~3.12(ii)]{BollobasJansonRiordan2007}
and the weighted tail condition.
The kernel bound gives $1-\rho_i\geq e^{-\kappa_*}$, and
$\mu\rho,\mu(1-\rho)\in\mathcal H_w$ because both are
coordinatewise bounded by $\mu$.
The empirical tail condition also gives
$\sup_n\sum_iw_i\mu_i^{(n)}<\infty$.
Cauchy--Schwarz gives, for $x\in\mathcal H_w$,
\begin{equation}
 \|x\|_1\leq C_w^{1/2}\|x\|_w,
 \qquad
 \|\kappa x\|_\infty\leq\kappa_*C_w^{1/2}\|x\|_w.
 \label{eq:Hw-embeddings}
\end{equation}
Hence $\mu^{(n)}\to\mu$ in $\ell^1$. By symmetry,
\begin{equation}
 \sup_j|\Delta_{n,ij}|
 =\sup_j|\Delta_{n,ji}|
 \leq\sqrt{w_i}\,\|\Delta_n\|_{w,\infty}\longrightarrow0
 \qquad(i\hbox{ fixed}).
 \label{eq:countable-row-convergence}
\end{equation}
In particular, the kernels converge at every pair of types. Since the
type space is discrete, this verifies the continuity requirement for a
graphical sequence in \cite[Definition~2.9]{BollobasJansonRiordan2007}.
Uniform boundedness and $\ell^1$ convergence of the empirical measures give
\begin{align*}
 &\left|\sum_{i,j}\kappa_{n,ij}\mu_i^{(n)}\mu_j^{(n)}
             -\sum_{i,j}\kappa_{ij}\mu_i\mu_j\right|\\
 &\qquad\leq2\kappa_*\|\mu^{(n)}-\mu\|_1
       +\left|\sum_{i,j}\Delta_{n,ij}\mu_i\mu_j\right|
       \longrightarrow0.
\end{align*}
The last term tends to zero by dominated convergence. Writing $e(G_n)$ for the number of edges, the exponential edge
probabilities satisfy
\[
 \frac1n\E e(G_n)
 =\frac12\sum_{i,j}\kappa_{n,ij}\mu_i^{(n)}\mu_j^{(n)}+O(n^{-1})
 \longrightarrow\frac12\sum_{i,j}\kappa_{ij}\mu_i\mu_j.
\]
The error includes the omission of diagonal vertex pairs. Thus the graph
sequence is graphical; the exponential parametrization is covered by
\cite[Remark~2.4]{BollobasJansonRiordan2007}.
The cited theorems give
$L_{n,i}/n\to \mu_i\rho_i$ in probability for every fixed $i$ and
$|C_2^{(n)}|=\Op(\log n)$.  Since
$0\leq L_{n,i}/n\leq\mu_i^{(n)}$ and $0\leq \mu_i\rho_i\leq\mu_i$, for every
fixed $m$,
\[
 \left\|\frac{L_n}{n}-\mu\rho\right\|_w^2
 \leq\sum_{i\in\mathcal I_m}w_i
             \left(\frac{L_{n,i}}n-\mu_i\rho_i\right)^2
       +\sum_{i>m}w_i(\mu_i^{(n)}+\mu_i).
\]
The finite sum tends to zero in probability.  The tail tends to zero
uniformly in $n$ as $m\to\infty$ by the weighted tail condition and
$\sum_iw_i\mu_i<\infty$.  Thus $L_n/n\to \mu\rho$ in $\mathcal H_w$.
The embedding in \eqref{eq:Hw-embeddings} gives the $\ell^1$ convergence,
which also implies $|C_1^{(n)}|/n\to\bar\rho$ in probability.

\end{proof}

Under the identification of type measures with their mass vectors, the maps in
\eqref{eq:exploration-residuals} have the coordinate form
\begin{equation}
 h_n(x)=x+\mu^{(n)}\odot e^{-\kappa_nx}-\mu^{(n)},
 \qquad h(x)=x+\mu\odot e^{-\kappa x}-\mu,
 \qquad x\in\mathcal H_w.
 \label{eq:countable-residual-maps}
\end{equation}
The identity \eqref{eq:exploration-equation} applies to these coordinate
expressions.  The next lemma controls the derivative used to solve
\eqref{eq:exploration-equation} and
the first order term caused by the parameter perturbations.

\begin{lemma}[Differentiability and first order expansion]
\label{lem:countable-terminal-maps}
The maps in \eqref{eq:countable-residual-maps} take values in
$\mathcal H_w$, are continuously Fr\'echet differentiable, and satisfy $h(\mu\rho)=0$.
If $\mu^{(n)}\to\mu$ in $\mathcal H_w$ and
$\|\Delta_n\|_{w,\infty}\to0$, then, for every fixed $\delta>0$,
\begin{equation}
 \sup_{\|x-\mu\rho\|_w\leq\delta}\|Dh_n(x)-Dh(x)\|_{\mathrm{op}}
 \leq C_\delta\bigl(\|\mu^{(n)}-\mu\|_w
                   +\|\Delta_n\|_{w,\infty}\bigr)\longrightarrow0.
 \label{eq:countable-derivative}
\end{equation}
If, in addition, $\sqrt n(\mu^{(n)}-\mu)\to\beta$ in $\mathcal H_w$
and \eqref{eq:main-countable-kernels} holds, then
\begin{equation}
 \sqrt n\,h_n(\mu\rho)\longrightarrow-b,\qquad
 b:=D_\rho\beta+D_{\mu(1-\rho)}\Lambda(\mu\rho)
 \quad\hbox{in }\mathcal H_w.
 \label{eq:countable-drift}
\end{equation}
\end{lemma}

\begin{proof}
The fixed point equation \eqref{eq:countable-fixed-point} and
\eqref{eq:countable-residual-maps} give $h(\mu\rho)=0$.
To differentiate the maps, we first bound the kernel action and its
composition with coordinatewise multiplication.  The kernel norm gives
\[
 \|Bx\|_\infty\leq\|B\|_{w,\infty}\|x\|_w.
\]
For $a\in\mathcal H_w$,
\begin{equation}
 \|D_aB\|_{\mathrm{HS}}^2
 =\sum_iw_i a_i^2\sum_j\frac{B_{ij}^2}{w_j}
 \leq\|a\|_w^2\|B\|_{w,\infty}^2.
 \label{eq:countable-kernel-multiplier-bound}
\end{equation}
In particular, $D_aB$ is a bounded linear operator whenever
$a\in\mathcal H_w$ and $\|B\|_{w,\infty}<\infty$.
For $\|x\|_w\leq R$ and $\|u\|_w\leq1$, the scalar exponential remainder
gives, uniformly in $n$,
\[
 \|h_n(x+u)-h_n(x)-u+
       D_{\mu^{(n)}\odot e^{-\kappa_nx}}\kappa_nu\|_w
 \leq C_R\|\mu^{(n)}\|_w\|u\|_w^2.
\]
Indeed $\|\kappa_nu\|_\infty\leq C\|u\|_w$, and the exponential
and its first two derivatives are bounded on the resulting fixed scalar
intervals.  The same bound holds for $h$ with $\mu$ in place of
$\mu^{(n)}$.  Thus both maps are Fr\'echet differentiable, with
\begin{equation}
 Dh_n(x)=I-D_{\mu^{(n)}\odot\exp\{-\kappa_nx\}}\kappa_n,
 \qquad
 Dh(x)=I-D_{\mu\odot\exp\{-\kappa x\}}\kappa.
 \label{eq:countable-map-derivatives}
\end{equation}
The operator bound \eqref{eq:countable-kernel-multiplier-bound} and local
Lipschitz continuity of the exponential show that the derivatives in
\eqref{eq:countable-map-derivatives} are continuous.
Under the parameter convergence assumptions, put
$F_n(x)=e^{-\kappa_nx}$ and $F(x)=e^{-\kappa x}$, the same exponential
factors used for the measure maps in
Lemma~\ref{lem:compact-terminal-maps}.  On the ball $\|x-\mu\rho\|_w\leq\delta$,
for fixed $\delta>0$,
\begin{align*}
 \|\mu^{(n)}\odot F_n(x)-\mu\odot F(x)\|_w
 &\leq C_\delta\bigl(\|\mu^{(n)}-\mu\|_w+\|\Delta_n\|_{w,\infty}\bigr),\\
 Dh_n(x)-Dh(x)
 &=-D_{\mu^{(n)}\odot F_n(x)-\mu\odot F(x)}\kappa
   -D_{\mu^{(n)}\odot F_n(x)}\Delta_n.
\end{align*}
Taking operator norms, using
$\|D_aB\|_{\mathrm{op}}\leq\|a\|_w\|B\|_{w,\infty}$, gives
\eqref{eq:countable-derivative}.

It remains to calculate $\sqrt n\,h_n(\mu\rho)$.
Under \eqref{eq:main-countable-counts} and
\eqref{eq:main-countable-kernels}, put $\delta_n=\Delta_n(\mu\rho)$.
The kernel expansion \eqref{eq:main-countable-kernels} gives
\[
 \|\delta_n\|_\infty=O(n^{-1/2}),\qquad
 \|\sqrt n\,\delta_n-\Lambda(\mu\rho)\|_\infty
 \leq\sqrt n\,\|E^{(n)}\|_{w,\infty}\|\mu\rho\|_w\longrightarrow0.
\]
Using $e^{-\kappa(\mu\rho)}=(1-\rho)$, which follows from
\eqref{eq:countable-fixed-point}, expand $e^{-\delta_n}=1-\delta_n+O(\|\delta_n\|_\infty^2)$
coordinatewise.  Since $h(\mu\rho)=0$, this gives
\begin{align*}
 h_n(\mu\rho)&=-D_\rho(\mu^{(n)}-\mu)-D_{\mu(1-\rho)}\delta_n+\epsilon_n,\\
 \|\epsilon_n\|_w
 &\leq C\bigl(\|\mu\|_w\|\delta_n\|_\infty^2
       +\|\mu^{(n)}-\mu\|_w\|\delta_n\|_\infty\bigr)=O(n^{-1}).
\end{align*}
Consequently, \eqref{eq:countable-drift} holds.

\end{proof}

To solve \eqref{eq:exploration-equation} near $\mu\rho$, we need a bounded
inverse of $Dh(\mu\rho)$.  The remaining type masses $\mu_i(1-\rho_i)$ determine a
symmetric operator on $\ell^2$:
\[
 A_{ij}=\sqrt{\mu_i(1-\rho_i)}\,\kappa_{ij}\sqrt{\mu_j(1-\rho_j)},\qquad
 r=\|A\|_{\ell^2\to\ell^2}.
\]
The next lemma proves $r<1$ and gives a quantitative inverse bound.

\begin{lemma}[Invertibility of the limiting derivative]
\label{lem:countable-inverse}
We have $r<1$, and
\begin{equation}
 J=Dh(\mu\rho)=I-D_{\mu(1-\rho)}\kappa
 \label{eq:countable-invertibility}
\end{equation}
is boundedly invertible on $\mathcal H_w$, with
\begin{equation}
 \|J^{-1}\|_{\mathrm{op}}
 \leq 1+
 \frac{\kappa_* C_w^{1/2}}{1-r}
 \left(\sum_i\mu_i(1-\rho_i)\right)^{1/2}
 \left(\sup_i w_i\mu_i(1-\rho_i)\right)^{1/2}.
 \label{eq:countable-inverse-bound}
\end{equation}
\end{lemma}

\begin{proof}
\noindent Step 1: A positive comparison vector.
The kernel bound and \eqref{eq:countable-fixed-point} give
$0\leq\rho_i\leq1-e^{-\kappa_*}<1$.
Supercriticality makes the survival vector nonzero.  If
$S=\{i:\rho_i>0\}$ were a proper subset, the fixed point equation
would give $\kappa_{ij}=0$ for $i\notin S$, $j\in S$; since
$\mu_i>0$ for every $i$, this would contradict irreducibility.
Thus $0<\rho_i<1$ for every $i$.

The estimate
$\sum_{i,j}A_{ij}^2\leq\kappa_*^2(\sum_i \mu_i(1-\rho_i))^2<\infty$
makes $A$ selfadjoint and compact.  Define the positive comparison vector
\[
 t_i=\frac{\mu_i\rho_i}{\sqrt{\mu_i(1-\rho_i)}}
     =\frac{\sqrt{\mu_i}\rho_i}{\sqrt{1-\rho_i}}.
\]
It belongs to $\ell^2$, because $1-\rho_i\geq e^{-\kappa_*}$.
The fixed point equation \eqref{eq:countable-fixed-point} gives
\[
 (At)_i=\sqrt{\mu_i(1-\rho_i)}\,[-\log(1-\rho_i)]<t_i,
\]
where the strict inequality follows from $-(1-\rho_i)\log(1-\rho_i)<\rho_i$.
\par\smallskip\noindent Step 2: The spectral gap.
Since $A$ has nonnegative matrix entries,
$|\langle f,Af\rangle_{\ell^2}|
 \leq\langle|f|,A|f|\rangle_{\ell^2}$.
Compact selfadjointness therefore gives a nonzero $g\geq0$ such that
$Ag=rg$, where $r=\|A\|>0$.  Here $r>0$ because supercriticality
makes the kernel nonzero and $1-\rho_i\geq e^{-\kappa_*}$.
Consequently,
\[
 (1-r)\langle g,t\rangle_{\ell^2}
 =\langle g,(I-A)t\rangle_{\ell^2}>0,
\]
and hence $r<1$.

\par\smallskip\noindent Step 3: The inverse on the weighted space.
To carry this inverse to $\mathcal H_w$, define
\[
 \begin{gathered}
 U:\ell^2\to\mathcal H_w,\qquad (Uf)_i=\sqrt{\mu_i(1-\rho_i)}f_i,\\
 V:\mathcal H_w\to\ell^2,\qquad (Vx)_i=\sqrt{\mu_i(1-\rho_i)}(\kappa x)_i.
 \end{gathered}
\]
These operators are bounded, since
\[
 \begin{gathered}
 \|U\|_{\ell^2\to\mathcal H_w}^2=\sup_i w_i\mu_i(1-\rho_i)
     \leq\sum_iw_i\mu_i<\infty,\\
 \|V\|_{\mathcal H_w\to\ell^2}
     \leq\kappa_*\left(C_w\sum_i\mu_i(1-\rho_i)\right)^{1/2}
 \end{gathered}
\]
by \eqref{eq:Hw-embeddings}.  Their products satisfy
$VU=A$ and $UV=D_{\mu(1-\rho)}\kappa$, so $J=I-UV$.
Direct multiplication now verifies
\[
 J^{-1}=I+U(I-A)^{-1}V.
\]
The Neumann series gives $\|(I-A)^{-1}\|\leq(1-r)^{-1}$,
which proves \eqref{eq:countable-inverse-bound}.
\end{proof}

\subsection{The limit for deterministic types}\label{sec:countable-martingale}

Lemmas~\ref{lem:countable-seed-approximation}--%
\ref{lem:countable-terminal-equation} approximate the giant by the
explored vector and identify the noise in
\eqref{eq:exploration-equation}. Their hypotheses are the ordinary
parameter convergence conditions of Lemma~\ref{lem:countable-components},
which also hold after conditioning on independently sampled types.

We use the seeded exploration of Subsection~\ref{sec:exploration},
with $0\le\kappa_n\le\kappa_*$.  For countable types, indicator functions
give its coordinate martingales.  Their finite dimensional limit is the first step toward the
full weighted Hilbert space limit; uniform coordinate variance bounds
will later control the omitted types.

For a finite or countable type set, let
$V_i^{(n)}=\{v\in[n]:x_v^{(n)}=i\}$ and identify the measures in \eqref{eq:exploration-measures} with their
coordinate sequences.  Write
$r_{n,t,i}=|R_{n,t}\cap V_i^{(n)}|$, $u_{n,t,i}=|U_{n,t}\cap V_i^{(n)}|$, and
$s_{n,i}=|I_{n,0}\cap V_i^{(n)}|$.  Put $r_{n,t}=(r_{n,t,i})_i$,
$u_{n,t}=(u_{n,t,i})_i$, and $s_n=(s_{n,i})_i$, and write
\begin{equation}
 Q_{n,t,i}=e^{-(\kappa_nr_{n,t}/n)_i},\qquad
 W_{n,t,i}=u_{n,t,i}/Q_{n,t,i},\qquad
 M_{n,t,i}=W_{n,t,i}/n.
 \label{eq:countable-exploration-coordinates}
\end{equation}
Write $W_{n,\cdot,i}=(W_{n,t,i})_{0\leq t\leq n}$.  The predictable
quadratic variation of this coordinate martingale is denoted by
\[
 \langle W_{n,\cdot,i}\rangle_t
 =\langle W_{n,\cdot}(\ind_{\{i\}}),W_{n,\cdot}(\ind_{\{i\}})\rangle_t.
\]

\begin{lemma}[Limits of finitely many type coordinates]
\label{lem:martingale-clt}
Fix finitely many type indices $i=1,\ldots,m$ and constants
$\mu_i\geq0$, $\rho_i\in[0,1)$.  Suppose that
$\mu_i^{(n)}=n_{i}^{(n)}/n\xrightarrow{\mathrm p}\mu_i$ and
$Q_{n,\tau_n,i}\xrightarrow{\mathrm p}1-\rho_i$ for these indices.  Then
\begin{equation}
 \left(\frac{W_{n,\tau_n,i}-W_{n,0,i}}{\sqrt n}:1\leq i\leq m\right)
 \Rightarrow\Normal_m\left(0,
       \diag(\frac{\mu_i\rho_i}{1-\rho_i}:1\leq i\leq m)\right).
 \label{eq:finite-coordinate-martingale-limit}
\end{equation}
Moreover, uniformly in the type index,
\begin{equation}
 \E(\langle W_{n,\cdot,i}\rangle_{\tau_n}\mid\mathcal E_n)\leq C n_i^{(n)}.
 \label{eq:coordinate-bracket-bound}
\end{equation}
\end{lemma}

\begin{proof}
We first bound the coordinate martingales to replace $M_{n,t,i}$ by its
limiting type proportion.  The bracket then telescopes, and the
Lindeberg bound \eqref{eq:common-lindeberg} completes the central limit
theorem.

At termination $I_{n,\tau_n}=\varnothing$, so the type counts satisfy
\begin{equation}
 r_{n,\tau_n}+u_{n,\tau_n}=(n_{i}^{(n)})_i.
 \label{eq:conservation}
\end{equation}
For $t<\tau_n$, put $B_{n,t,i}=|B_{n,t}\cap V_i^{(n)}|$.
The conditional independence in Lemma~\ref{lem:exploration-martingales}
shows that, given $\mathcal F_{n,t}$, the variables $B_{n,t,i}$ are independent
with laws $\operatorname{Bin}(u_{n,t,i},p_{n,t}(i))$.
Taking $f=\ind_{\{i\}}$ in \eqref{eq:measure-increment} gives
\begin{equation}
 \Delta W_{n,t,i}
 =\frac{u_{n,t,i}p_{n,t}(i)-B_{n,t,i}}{Q_{n,t,i}(1-p_{n,t}(i))}.
 \label{eq:increment}
\end{equation}
Consequently,
\begin{align}
 \E((\Delta W_{n,t,i})^2\mid\mathcal F_{n,t})
 &=\frac{u_{n,t,i}p_{n,t}(i)}{Q_{n,t,i}^2(1-p_{n,t}(i))},
 \label{eq:conditional-variance}\\
 \E(\Delta W_{n,t,i}\Delta W_{n,t,j}\mid\mathcal F_{n,t})&=0
 \quad(i\ne j).
 \label{eq:cross-variation}
\end{align}
The conditional variance in \eqref{eq:conditional-variance} is at most
$C n_i^{(n)}/n$, since
$p_{n,t}(i)\leq\kappa_*/n$, $u_{n,t,i}\leq n_{i}^{(n)}$, and
$Q_{n,t,i},1-p_{n,t}(i)\geq e^{-\kappa_*}$.
Summing over $\tau_n\leq n$ proves \eqref{eq:coordinate-bracket-bound}.
The martingale isometry and Doob's inequality, conditional on
$\mathcal E_n\subset\mathcal F_{n,0}$, yield
\begin{align*}
 \E\left[\sup_{t\leq\tau_n}|M_{n,t,i}-M_{n,0,i}|^2
          \biggm|\mathcal E_n\right]
 &\leq\frac4{n^2}
       \E[(W_{n,\tau_n,i}-W_{n,0,i})^2\mid\mathcal E_n]\\
 &=\frac4{n^2}\E[\langle W_{n,\cdot,i}\rangle_{\tau_n}\mid\mathcal E_n]
 \leq C n_i^{(n)}/n^2.
\end{align*}
Because $M_{n,0,i}=\mu_i^{(n)}-s_{n,i}/n$ and
$s_{n,i}\leq\ell_n$, it implies
$\sup_{t\leq\tau_n}|M_{n,t,i}-\mu_i|\to0$ in probability for every fixed $i$.
Equation \eqref{eq:conditional-variance} now yields
\begin{align*}
 \frac1n\langle W_{n,\cdot,i}\rangle_{\tau_n}
 &=\sum_{t<\tau_n}M_{n,t,i}
       (Q_{n,t+1,i}^{-1}-Q_{n,t,i}^{-1})\\
 &=\mu_i(Q_{n,\tau_n,i}^{-1}-1)+\op(1)
 \xrightarrow{\mathrm p}\frac{\mu_i\rho_i}{1-\rho_i}.
\end{align*}
The error from replacing $M_{n,t,i}$ by $\mu_i$ is at most
\[
 (e^{\kappa_*}-1)\sup_{t\leq\tau_n}|M_{n,t,i}-\mu_i|
 \xrightarrow{\mathrm p}0.
\]
For every $a=(a_1,\ldots,a_m)^{\mathsf T}\in\R^m$, the vanishing
cross covariations in \eqref{eq:cross-variation} give
\[
 \frac1n\left\langle\sum_{i=1}^m a_iW_{n,\cdot,i}\right\rangle_{\tau_n}
 \xrightarrow{\mathrm p}\sum_{i=1}^m a_i^2\frac{\mu_i\rho_i}{1-\rho_i}.
\]
Apply \eqref{eq:common-lindeberg} to
$f=\sum_{i=1}^m a_i\ind_{\{i\}}$.  This proves conditional Lindeberg for
$a^{\mathsf T}(W_{n,t,i}:1\leq i\leq m)/\sqrt n$.
The scalar martingale array central limit theorem
\cite[Corollary~3.1 and the following remark, pp.~58--59]{HallHeyde1980},
followed by the Cram\'er--Wold device, proves the asserted limit.
\end{proof}

Throughout the next three lemmas, assume the hypotheses of
Lemma~\ref{lem:countable-components} and run the seeded exploration of
Subsection~\ref{sec:exploration}.

\begin{lemma}[Approximation by the seeded exploration]
\label{lem:countable-seed-approximation}
The seeds meet the giant with probability tending to one, and
\begin{equation}
 \frac{\|r_{n,\tau_n}-L_n\|_w}{\sqrt n}
 \xrightarrow{\mathrm p}0,
 \qquad
 \frac{r_{n,\tau_n}}n\xrightarrow{\mathrm p}\mu\rho
 \quad\hbox{in both }\mathcal H_w\hbox{ and }\ell^1.
 \label{eq:countable-seeded-approximation}
\end{equation}
\end{lemma}

\begin{proof}
By \eqref{eq:countable-first-order}, $|C_1^{(n)}|/n\to\bar\rho>0$
in probability.  Applying \eqref{eq:seed-hit} with $a=\bar\rho/2$
shows that the seeds meet the giant with probability tending to one.
We next control the additional vertices in the weighted norm.

Let $D_n$ be the union of the components other than $C_1^{(n)}$ that
intersect $I_{n,0}$, and put
$d_{n,i}=|D_n\cap V_i^{(n)}|$ and $d_n=(d_{n,i})_i$.
On the event that the seeds meet the giant,
$r_{n,\tau_n}=L_n+d_n$.  For every graph and seed realization,
\[
 |D_n|\leq\ell_n|C_2^{(n)}|,\qquad
 \|d_n\|_w^2\leq |D_n|\sum_{v\in D_n}w_{x_v^{(n)}}.
\]
A nongiant component $C$ contains a seed with conditional probability at
most $\ell_n|C|/n$, by the union bound and the inclusion probability
$\ell_n/n$ of each vertex.  Set
$m_{n,w}=\sum_iw_i\mu_i^{(n)}$, which is uniformly bounded by the
weighted tail condition and coordinate convergence.  Therefore
\begin{align*}
 \E\left[\sum_{v\in D_n}w_{x_v^{(n)}}\biggm|G_n\right]
 &\leq\frac{\ell_n|C_2^{(n)}|}{n}
             \sum_{v\notin C_1^{(n)}}w_{x_v^{(n)}}
 \leq\ell_n|C_2^{(n)}|m_{n,w},\\
 \E(\|d_n\|_w^2\mid G_n)
 &\leq\ell_n^2|C_2^{(n)}|^2m_{n,w}.
\end{align*}
For every $a,R>0$, conditional Markov's inequality gives
\[
 \Pp(\|d_n\|_w>a\sqrt n)
 \leq\Pp(|C_2^{(n)}|>R\log n)
       +\frac{CR^2(\log n)^4}{a^2n}.
\]
First let $n\to\infty$ and then $R\to\infty$.
Since $|C_2^{(n)}|=\Op(\log n)$ by
Lemma~\ref{lem:countable-components}, it follows that
$\|d_n\|_w/\sqrt n\to0$ in probability.  Hence
\[
 \Pp(\|r_{n,\tau_n}-L_n\|_w>a\sqrt n)
 \leq\Pp(I_{n,0}\cap C_1^{(n)}=\varnothing)
       +\Pp(\|d_n\|_w>a\sqrt n)\longrightarrow0.
\]
On the event that the seeds meet the giant, also
$\|r_{n,\tau_n}-L_n\|_1=|D_n|\leq\ell_n|C_2^{(n)}|$.
The two asserted laws of large numbers now follow from
Lemma~\ref{lem:countable-components} and
\eqref{eq:Hw-embeddings}.
\end{proof}

The fluctuations of $W_{n,t}$ around its initial value carry the edge
noise.  The coordinate projection $P_m$ from
\eqref{eq:countable-projection} separates finitely many coordinates of
its normalized increment $N_n$ from the remaining tail.
The finite coordinate martingale variances identify the candidate
Gaussian law.  The proof below also controls the omitted coordinates,
which is needed for convergence in $\mathcal H_w$.

\begin{lemma}[The Hilbert space martingale limit]
\label{lem:countable-martingale-limit}
The normalized martingale satisfies
\begin{equation}
 \sup_{t\leq\tau_n}\|M_{n,t}-\mu^{(n)}\|_w
 =\Op(n^{-1/2}).
 \label{eq:countable-M}
\end{equation}
Moreover, $Q_{n,\tau_n,i}\to1-\rho_i$ in probability for each fixed $i$, and
\[
 N_n\Rightarrow\Normal\!\left(0,D_{w\mu\rho/(1-\rho)}\right)
 \quad\hbox{in }\mathcal H_w.
\]
\end{lemma}

\begin{proof}
We first show that the normalized martingale remains close to its initial
type proportions.  This will identify its quadratic variation at $\tau_n$.
For each $n$, the processes $M_{n,t}$ and $W_{n,t}$ are supported on the
finite set of types present in $G_n$.  Lemma~\ref{lem:exploration-martingales}
therefore makes them square integrable $\mathcal H_w$ valued
martingales.  Also $M_{n,0}=\mu^{(n)}-s_n/n$ and
$e^{-\kappa_*}\leq Q_{n,t,i}\leq1$.
Summing the coordinate conditional variances in
\eqref{eq:conditional-variance} gives
\begin{align*}
 \sum_{t<n}\E\bigl(\|\Delta M_{n,t}\|_w^2
                          \mid\mathcal F_{n,t}\bigr)
 &\leq\frac Cn\sum_iw_i\mu_i^{(n)}=O(n^{-1}),\\
 \E\left\|\frac{s_n}{n}\right\|_w^2
 &\leq\frac{\ell_n}{n^2}\E\sum_iw_is_{n,i}
 \leq\frac{C\ell_n^2}{n^2}.
\end{align*}
For the seed estimate, use
$\|s_n\|_w^2\leq\ell_n\sum_iw_is_{n,i}$ and
$\E s_{n,i}=\ell_n\mu_i^{(n)}$.
The Hilbert space martingale isometry and Doob's inequality yield
\[
 \E\sup_{t\leq\tau_n}\|M_{n,t}-M_{n,0}\|_w^2\leq C/n.
\]
Together with the seed estimate, this gives
$\E\sup_{t\leq\tau_n}\|M_{n,t}-\mu^{(n)}\|_w^2\leq C/n$,
which proves the rate in \eqref{eq:countable-M}.

For every fixed $i$, Lemma~\ref{lem:countable-seed-approximation} and
\eqref{eq:countable-row-convergence} give
\begin{align*}
 \left|(\kappa_nr_{n,\tau_n}/n)_i-(\kappa(\mu\rho))_i\right|
 &\leq\sup_j|\kappa_{n,ij}-\kappa_{ij}|
          \|r_{n,\tau_n}/n\|_1
       +\kappa_*\|r_{n,\tau_n}/n-\mu\rho\|_1
 \xrightarrow{\mathrm p}0.
\end{align*}
The exponential formula \eqref{eq:countable-exploration-coordinates} for
$Q_{n,\tau_n,i}$ and the fixed point equation
\eqref{eq:countable-fixed-point} imply $Q_{n,\tau_n,i}\to1-\rho_i$ in probability.
Applying \eqref{eq:finite-coordinate-martingale-limit} to the indices in
$\mathcal I_m$ gives
\[
 \left(N_{n,i}:i\in\mathcal I_m\right)
 \Rightarrow
 \Normal_{|\mathcal I_m|}\!\left(0,
       \diag(\frac{\mu_i\rho_i}{1-\rho_i}:i\in\mathcal I_m)\right).
\]
The diagonal covariance operator is positive and trace class, since
\begin{equation}
 \operatorname{tr}D_{w\mu\rho/(1-\rho)}
 =\sum_iw_i\frac{\mu_i\rho_i}{1-\rho_i}
 \leq(e^{\kappa_*}-1)\sum_iw_i\mu_i<\infty.
 \label{eq:countable-gaussian-moment}
\end{equation}
It therefore defines a Gaussian law on $\mathcal H_w$.

To pass to the whole space, the bound in
Lemma~\ref{lem:martingale-clt} and the martingale isometry give
\[
 \E N_{n,i}^2
 =\frac1n\E\langle W_{n,\cdot,i}\rangle_{\tau_n}
 \leq C\mu_i^{(n)}.
\]
Consequently,
\begin{equation}
 \sup_n\E\left\|(I-P_m)
                 N_n\right\|_w^2
 \leq C\sup_n\sum_{i>m}w_i\mu_i^{(n)}\longrightarrow0.
 \label{eq:countable-noise-tail}
\end{equation}
The finite coordinate limit and \eqref{eq:countable-noise-tail}
verify Lemma~\ref{lem:hilbert-projections} for this Gaussian law; its
coordinate tails have second moments tending to zero by
\eqref{eq:countable-gaussian-moment}.  Hence
\begin{equation}
 N_n\Rightarrow\Normal\!\left(0,D_{w\mu\rho/(1-\rho)}\right)
 \quad\hbox{in }\mathcal H_w.
 \label{eq:countable-W-limit}
\end{equation}
\end{proof}

Undoing the compensation in $N_n$ returns to the actual unexplored
counts.  In coordinates, the noise from \eqref{eq:exploration-noise} is
\begin{equation}
 Z_{n,i}=\frac{u_{n,\tau_n,i}-n\mu_i^{(n)}Q_{n,\tau_n,i}}{\sqrt n}
 =Q_{n,\tau_n,i}N_{n,i}
  -Q_{n,\tau_n,i}\frac{s_{n,i}}{\sqrt n}.
 \label{eq:countable-exploration-noise}
\end{equation}
This is the coordinate representation of \eqref{eq:exploration-noise}:
$Z_n$ measures the deviation of the remaining type counts from
$n\mu^{(n)}\odot Q_{n,\tau_n}$ on the $\sqrt n$ scale, and
\eqref{eq:exploration-equation} relates it to the explored type vector.
\begin{lemma}[Noise in the exploration equation]
\label{lem:countable-terminal-equation}
The exploration noise satisfies
\begin{equation}
 Z_n=D_{1-\rho}N_n+\op(1)
 \Rightarrow\Normal\!\left(0,D_{w\mu\rho(1-\rho)}\right)
 \quad\hbox{in }\mathcal H_w.
 \label{eq:countable-terminal-noise}
\end{equation}
The limiting coordinate variances are $\mu_i\rho_i(1-\rho_i)$, and
its coordinates are independent.
\end{lemma}

\begin{proof}
The seed term vanishes because $0<Q_{n,\tau_n,i}\leq1$ and
\[
 \E\left\|\frac{s_n}{\sqrt n}\right\|_w^2
 \leq C\ell_n^2/n\longrightarrow0.
\]
For each fixed $m$, Lemma~\ref{lem:countable-martingale-limit} gives
$P_m(D_{Q_{n,\tau_n}}-D_{1-\rho})N_n\to0$ in probability.
Indeed the finitely many multipliers converge to zero, and the finite
dimensional vectors $P_mN_n$ are tight.  Both diagonal operators are
contractions and commute with $P_m$, so
\[
 \E\|(I-P_m)(D_{Q_{n,\tau_n}}-D_{1-\rho})N_n\|_w^2
 \leq4\E\|(I-P_m)N_n\|_w^2.
\]
By \eqref{eq:countable-noise-tail}, the right side vanishes uniformly
in $n$ as $m\to\infty$.  Thus
$(D_{Q_{n,\tau_n}}-D_{1-\rho})N_n\to0$ in probability in
$\mathcal H_w$.  The vanishing seed term gives the stated
approximation by $D_{1-\rho}N_n$.  The bounded diagonal map $D_{1-\rho}$
transforms the covariance in \eqref{eq:countable-W-limit} into
\[
 D_{1-\rho}D_{w\mu\rho/(1-\rho)}D_{1-\rho}
 =D_{w\mu\rho(1-\rho)},
\]
which proves \eqref{eq:countable-terminal-noise}.
\end{proof}

\begin{proof}[Proof of Theorem~\ref{thm:countable-concrete}]
The exploration equation determines the fluctuation once its random error
and its response to a change in the explored mass are known.  We first
locate the explored mass near $\mu\rho$, then linearize that equation and
transfer the result to the giant.

\noindent Step 1: Approximate the giant by the exploration.
The parameter assumptions give $\mu^{(n)}\to\mu$ in $\mathcal H_w$ and
\[
 \|\Delta_n\|_{w,\infty}
 \leq n^{-1/2}\|\Lambda\|_{w,\infty}
       +\|E^{(n)}\|_{w,\infty}\longrightarrow0.
\]
Together with the assumed uniform tail bound, these are the hypotheses
of Lemma~\ref{lem:countable-components}.  Thus the giant has limiting
density $\bar\rho>0$ and $|C_2^{(n)}|=\Op(\log n)$.
The seed miss probability in \eqref{eq:seed-hit} tends to zero with
$a=\bar\rho/2$: its first term vanishes by the giant law of large numbers,
and $(1-\bar\rho/2)^{\ell_n}\to0$ because $\ell_n\to\infty$.
On the event that the seeds meet the giant, the exploration contains the
giant and only additional components meeting the seeds.
Lemma~\ref{lem:countable-seed-approximation} controls these additional
vertices in the weighted norm and gives
\[
 \widehat\mu_n\xrightarrow{\mathrm p}\mu\rho,
 \qquad
 \frac{r_{n,\tau_n}-L_n}{\sqrt n}\xrightarrow{\mathrm p}0
 \quad\hbox{in }\mathcal H_w.
\]
The first convergence selects the nonzero fixed point about which to
expand; the second permits the replacement of the giant by the explored
vector on the fluctuation scale.

\par\smallskip\noindent Step 2: Separate edge noise from parameter drift.
Conservation of type mass gives the exact relation
\eqref{eq:exploration-equation}.  Lemma~\ref{lem:countable-terminal-equation}
identifies its error:
\[
 \sqrt n\,h_n(\widehat\mu_n)=-Z_n,
 \qquad Z_n\Rightarrow\Normal\!\left(0,D_{w\mu\rho(1-\rho)}\right)
 \quad\hbox{in }\mathcal H_w.
\]
At the limiting mass $\mu\rho$, the
fixed point equation gives $h(\mu\rho)=0$, while the finite parameter
changes give, by \eqref{eq:countable-drift},
\[
 \sqrt n\,h_n(\mu\rho)\longrightarrow-b,
 \qquad b=D_\rho\beta+D_{\mu(1-\rho)}\Lambda(\mu\rho).
\]
The two summands of $b$ come from the perturbations of the type proportions
and of the kernel, respectively.  Thus the residual at the explored mass
is random on the $n^{-1/2}$ scale, and its value at the limiting mass has
a deterministic shift on the same scale.  In particular,
$\sqrt n\,h_n(\widehat\mu_n)=\Op(1)$ and
$\sqrt n\,h_n(\mu\rho)=O(1)$.

\par\smallskip\noindent Step 3: Linearize and solve for the explored mass.
Lemma~\ref{lem:countable-terminal-maps} shows that $h_n,h$ are continuously
Fr\'echet differentiable and gives local uniform convergence of their
derivatives in \eqref{eq:countable-derivative}.  By
Lemma~\ref{lem:countable-inverse}, their limiting derivative
$J=Dh(\mu\rho)$ has a bounded inverse.
Consequently, Lemma~\ref{lem:linearization} applies with
$g_n=h_n$, $g=h$, $x_\star=\mu\rho$, and $X_n=\widehat\mu_n$.
It first yields $\sqrt n(\widehat\mu_n-\mu\rho)=\Op(1)$ and then
\[
 \sqrt n\{h_n(\widehat\mu_n)-h_n(\mu\rho)\}
 =J\sqrt n(\widehat\mu_n-\mu\rho)+\op(1).
\]
Substituting the two residual values
from Step~2 gives
\[
 J\sqrt n(\widehat\mu_n-\mu\rho)=b-Z_n+\op(1),
\]
and boundedness of $J^{-1}$ therefore gives
\[
 \sqrt n(\widehat\mu_n-\mu\rho)
 =J^{-1}(b-Z_n)+\op(1)
 \quad\hbox{in }\mathcal H_w.
\]
Thus the parameter shift and the edge noise are transformed by the same
inverse derivative.

\par\smallskip\noindent Step 4: Return to the giant and identify its Gaussian law.
Since $r_{n,\tau_n}=n\widehat\mu_n$, Step~1 and Slutsky's theorem give
\[
 \frac{L_n-n\mu\rho}{\sqrt n}
 =\sqrt n(\widehat\mu_n-\mu\rho)
    +\frac{L_n-r_{n,\tau_n}}{\sqrt n}
 =J^{-1}(b-Z_n)+\op(1).
\]
The Gaussian law in Step~2 therefore gives limiting mean $J^{-1}b$.
The weighted inner product accounts for the factor $w$ in its noise
covariance: the coordinate variances $\mu_i\rho_i(1-\rho_i)$ give
\[
 \sum_iw_i^2\mu_i\rho_i(1-\rho_i)x_iy_i
 =\langle D_{w\mu\rho(1-\rho)}x,y\rangle_w,
 \qquad x,y\in\mathcal H_w.
\]
This diagonal operator is trace class, since its trace in the orthonormal
basis $(\mathbf e_i/\sqrt{w_i})_i$ is
\[
 \sum_iw_i\mu_i\rho_i(1-\rho_i)\leq\sum_iw_i\mu_i<\infty.
\]
Applying the bounded map $-J^{-1}$ preserves the trace-class property
and gives covariance
\[
 J^{-1}D_{w\mu\rho(1-\rho)}(J^{-1})^{*_{w}}.
\]
Consequently,
\[
 \frac{L_n-n\mu\rho}{\sqrt n}
 \Rightarrow\Normal\!\left(J^{-1}b,
 J^{-1}D_{w\mu\rho(1-\rho)}(J^{-1})^{*_{w}}\right).
\]
Substituting the expression for $b$ in \eqref{eq:countable-drift}
proves \eqref{eq:concrete-countable-limit}.
\end{proof}

\subsection{Independently sampled types}\label{sec:countable-iid-proof}

Random type allocation and random edges contribute two separate
fluctuations.  We first prove their joint convergence and asymptotic
independence, then insert both terms into
\eqref{eq:exploration-equation}.

Throughout this subsection assume the hypotheses of
Corollary~\ref{cor:countable-iid}.  We use the empirical fluctuation $H_n$
from \eqref{eq:countable-iid-empirical-fluctuation} and the coordinate
projections from \eqref{eq:countable-projection}.

\begin{lemma}[Empirical type fluctuations and weighted tails]
\label{lem:countable-iid-types}
The empirical type vector satisfies
\[
 H_n\Rightarrow\Normal\!\left(0,D_{w\mu}-\mu\otimes_w\mu\right)
 \quad\hbox{in }\mathcal H_w.
\]
Moreover,
$\sum_iw_i|\mu_i^{(n)}-\mu_i|\to0$ and
$\lim_{m\to\infty}\sup_n\sum_{i>m}w_i\mu_i^{(n)}=0$ almost surely.
\end{lemma}

\begin{proof}
Every fixed finite projection of $H_n$ satisfies the multivariate
central limit theorem.  For each $i,j\in\mathcal S$, independence
and centering give
\[
\begin{aligned}
 \E[H_{n,i}H_{n,j}]
 =\frac1n\sum_{v=1}^n
   \E\!\left[
   (\ind_{\{X_v=i\}}-\mu_i)
   (\ind_{\{X_v=j\}}-\mu_j)
   \right]
 =\mu_i\ind_{\{i=j\}}-\mu_i\mu_j.
\end{aligned}
\]
Consequently,
\begin{equation}
 \E\|(I-P_m)H_n\|_w^2
 =\sum_{i>m}w_i\mu_i(1-\mu_i)
 \leq\sum_{i>m}w_i\mu_i
 \xrightarrow[m\to\infty]{}0,
 \label{eq:countable-iid-type-tail}
\end{equation}
uniformly in $n$.

The covariance operator of $\mathbf e_{X_1}-\mu$ is
$D_{w\mu}-\mu\otimes_w\mu$.  It is positive and has finite trace
$\sum_iw_i\mu_i(1-\mu_i)$, so it determines a Gaussian law on
$\mathcal H_w$.  Lemma~\ref{lem:hilbert-projections}, together with
\eqref{eq:countable-iid-type-tail}, proves convergence to this law.

We next prove the almost sure assertions, which allow us to fix
the types when studying the edge noise.  Since
$\E[w_{X_1}]=\sum_iw_i\mu_i<\infty$, the scalar strong law applies
to $w_{X_v}$ and to each coordinate indicator.  As $\mathcal S$
is countable, on a single event of probability one,
\[
 \mu_i^{(n)}\longrightarrow\mu_i
 \quad\text{for every }i\in\mathcal S,
 \qquad
 \sum_iw_i\mu_i^{(n)}
 =\frac1n\sum_{v=1}^nw_{X_v}
 \longrightarrow\sum_iw_i\mu_i.
\]
Applying Lemma~\ref{lem:countable-weighted-moments} on this event
yields
\begin{equation}
 \sum_iw_i|\mu_i^{(n)}-\mu_i|\longrightarrow0,
 \qquad
 \lim_{m\to\infty}\sup_n\sum_{i>m}w_i\mu_i^{(n)}=0
 \quad\text{almost surely}.
 \label{eq:countable-iid-type-lln}
\end{equation}
Finally, since $w_i\geq1$ and $|\mu_i^{(n)}-\mu_i|\leq1$,
\[
 \|\mu^{(n)}-\mu\|_{\ell^1}
 \leq\sum_iw_i|\mu_i^{(n)}-\mu_i|,
 \qquad
 \|\mu^{(n)}-\mu\|_w^2
 \leq\sum_iw_i|\mu_i^{(n)}-\mu_i|.
\]
Thus $\mu^{(n)}\to\mu$ almost surely in both $\ell^1$ and
$\mathcal H_w$.
\end{proof}

Here the environment is $\mathcal E_n=\sigma(X_1,\ldots,X_n)$.
Run the same seeded exploration, with $N_n$ and $Z_n$ as in
\eqref{eq:exploration-martingale-noise} and
\eqref{eq:countable-exploration-noise}.
The conditional Gaussian edge law has a deterministic limit, which gives
asymptotic independence of the empirical type fluctuation and the edge noise.

\begin{lemma}[Joint limit of type and edge noises]
\label{lem:countable-iid-joint}
In $\mathcal H_w\times\mathcal H_w$,
\[
 (H_n,Z_n)\Rightarrow\Normal\!\left(0,
 \begin{pmatrix}
 D_{w\mu}-\mu\otimes_w\mu&0\\
 0&D_{w\mu\rho(1-\rho)}
 \end{pmatrix}\right).
\]
In particular, the two coordinates are asymptotically independent.
\end{lemma}

\begin{proof}
For any deterministic type sequence satisfying
\eqref{eq:countable-iid-type-lln}, the hypotheses of
Lemma~\ref{lem:countable-components} hold with the fixed kernel $\kappa$.
Lemmas~\ref{lem:countable-seed-approximation}--\ref{lem:countable-terminal-equation}
therefore give
$Z_n\Rightarrow\Normal(0,D_{w\mu\rho(1-\rho)})$.
Consequently, conditioning on a type realization in the probability one
event of Lemma~\ref{lem:countable-iid-types}, we obtain, for every
finitely supported real sequence $a$,
\begin{equation}
 \phi_n(a):=
 \E\left[\exp\left\{\mathrm i\sum_i a_iZ_{n,i}\right\}
                   \biggm|\mathcal E_n\right]
 \longrightarrow
 \exp\left\{-\tfrac12\sum_i a_i^2\mu_i\rho_i(1-\rho_i)\right\}
 \quad\hbox{almost surely},
 \label{eq:countable-iid-conditional-characteristic}
\end{equation}
For any other finitely supported real sequence $b$, conditioning gives
\[
 \E\exp\left\{\mathrm i\sum_i b_iH_{n,i}
                    +\mathrm i\sum_i a_iZ_{n,i}\right\}
 =\E\left[\exp\left\{\mathrm i\sum_i b_iH_{n,i}\right\}
                 \phi_n(a)\right].
\]
By \eqref{eq:countable-iid-conditional-characteristic} and
$|\phi_n(a)|\leq1$, bounded convergence and the finite dimensional
limit of $H_n$ give
\[
 \exp\!\left\{-\tfrac12\left[
 \sum_i\mu_i b_i^2-\left(\sum_i\mu_i b_i\right)^2
 +\sum_i a_i^2\mu_i\rho_i(1-\rho_i)\right]\right\}.
\]
This is the characteristic function of the stated Gaussian law on
finite coordinates, whose two blocks are independent.

To pass from finite coordinates to the whole Hilbert space, average the
conditional bracket bound \eqref{eq:coordinate-bracket-bound} over the
types, and use the same seed sampling calculation as in the proof of
Lemma~\ref{lem:countable-martingale-limit}.  Since $\E\mu_i^{(n)}=\mu_i$ and
$\E s_{n,i}=\ell_n\mu_i$,
\begin{align*}
 \E\|(I-P_m)N_n\|_w^2
 &\leq C\sum_{i>m}w_i\mu_i,\\
 \E\|(I-P_m)s_n/\sqrt n\|_w^2
 &\leq\frac{\ell_n^2}{n}\sum_{i>m}w_i\mu_i.
\end{align*}
The diagonal multiplier $D_{Q_{n,\tau_n}}$ is a contraction and commutes with
$P_m$.  Since $\sup_n\ell_n^2/n<\infty$, these bounds imply
\begin{equation}
 \sup_n\E\|(I-P_m)Z_n\|_w^2
 \leq C\sum_{i>m}w_i\mu_i\longrightarrow0.
 \label{eq:countable-iid-edge-tail}
\end{equation}
The joint finite coordinate limits and the tail bounds
\eqref{eq:countable-iid-type-tail} and~\eqref{eq:countable-iid-edge-tail}
satisfy Lemma~\ref{lem:hilbert-projections} on
$\mathcal H_w\times\mathcal H_w$, proving the asserted joint limit
and independence.

\end{proof}

\begin{proof}[Proof of Corollary~\ref{cor:countable-iid}]
Lemma~\ref{lem:countable-iid-joint} identifies the joint limit of the
type and edge noises.  It remains to insert them into
\eqref{eq:exploration-equation}.  With the kernel fixed, substituting
$\mu\rho$ into \eqref{eq:countable-residual-maps} gives
\begin{equation}
 h_n(\mu\rho)=-D_\rho(\mu^{(n)}-\mu),
 \qquad \sqrt n\,h_n(\mu\rho)=-D_\rho H_n.
 \label{eq:countable-iid-residual-at-z}
\end{equation}
Conditioning on the probability one event in
Lemma~\ref{lem:countable-iid-types},
Lemma~\ref{lem:countable-seed-approximation} gives
$r_{n,\tau_n}/n\to\mu\rho$ in probability in $\mathcal H_w$ and
$(r_{n,\tau_n}-L_n)/\sqrt n\to0$ in probability in that space, first
conditionally and then after averaging over the types.
The derivative estimate \eqref{eq:countable-derivative}
holds almost surely by \eqref{eq:countable-iid-type-lln}, while
\eqref{eq:countable-iid-residual-at-z} expresses the first order term
in terms of the empirical fluctuation
\eqref{eq:countable-iid-empirical-fluctuation}.
The exploration identity \eqref{eq:exploration-equation} gives
$\sqrt n\,h_n(r_{n,\tau_n}/n)=-Z_n$.
Lemmas~\ref{lem:countable-inverse} and~\ref{lem:linearization} yield
\[
 \frac{L_n-n\mu\rho}{\sqrt n}
 =J^{-1}(D_\rho H_n-Z_n)+\op(1).
\]
Apply the bounded linear map $(h,u)\mapsto D_\rho h-u$ to the
joint Gaussian limit of Lemma~\ref{lem:countable-iid-joint}.
Its two independent blocks give covariance operator
\begin{align*}
 &D_\rho(D_{w\mu}-\mu\otimes_w\mu)D_\rho
       +D_{w\mu\rho(1-\rho)}\\
 &\qquad=D_{w\mu\rho}-(\mu\rho)\otimes_w(\mu\rho),
\end{align*}
with trace
$\sum_iw_i\mu_i\rho_i(1-\mu_i\rho_i)\leq\sum_iw_i\mu_i<\infty$.
Applying $J^{-1}$ therefore gives
\[
 \frac{L_n-n\mu\rho}{\sqrt n}
 \Rightarrow\Normal\!\left(0,
 J^{-1}\{D_{w\mu\rho}-(\mu\rho)\otimes_w(\mu\rho)\}
 (J^{-1})^{*_{w}}\right),
\]
which proves \eqref{eq:countable-iid-limit} and its finite type form
\eqref{eq:finite-iid-limit}.
\end{proof}

\subsection{The infinite rank example}\label{sec:countable-example}

We check that rounding and truncating the geometric and polynomial type
distributions produces errors smaller than $n^{-1/2}$ in the weighted
norm.  The kernel assumptions then follow from its positive constant
part.

\begin{proof}[Verification of Example~\ref{ex:countable-bounded-perturbations}]
For both type distributions, the rounding rule gives
$0\leq\mu_i-\mu_i^{(n)}<n^{-1}$ for $2\leq i\leq m_n$ and
$\mu_i^{(n)}=0$ for $i>m_n$.  Mass conservation gives
\[
 |\mu_1^{(n)}-\mu_1|
 \leq m_n/n+\sum_{i>m_n}\mu_i.
\]
Thus the first coordinate, rounding, and truncation contribute at most
\[
 nw_1\left(m_n/n+\sum_{i>m_n}\mu_i\right)^2,
 \qquad \frac1n\sum_{i=2}^{m_n}w_i,
 \qquad n\sum_{i>m_n}w_i\mu_i^2
\]
respectively to $n\|\mu^{(n)}-\mu\|_w^2$.
Substituting the geometric and polynomial choices from
Example~\ref{ex:countable-bounded-perturbations} gives
\begin{center}
\small
\begin{tabular}{@{}lcc@{}}
\hline
Error source & Geometric types & Polynomial types\\
\hline
First coordinate & $O((\log n)^2/n)$ & $O(n^{2/p-1})$\\
Rounding $2\leq i\leq m_n$ & $O(n^{\alpha-1})$
                          & $O(n^{(a+1)/p-1})$\\
Truncation $i>m_n$ & $O(n^{\alpha-1})$
                   & $O(n^{(a+1)/p-1})$\\
\hline
\end{tabular}
\end{center}
For geometric types and $n\ge2$, these bounds use $m_n=\lfloor\log_2n\rfloor$
and the geometric sums
\[
 \frac1n\sum_{i=2}^{m_n}2^{\alpha i}
 +n\sum_{i>m_n}2^{-(2-\alpha)i}=O(n^{\alpha-1}).
\]
For polynomial types, $m_n\asymp n^{1/p}$ and
$\sum_{i>m_n}\mu_i\leq C m_n^{1-p}$ give the first bound, while
\[
 \frac1n\sum_{i=2}^{m_n}i^a
 +nZ_p^{-2}\sum_{i>m_n}i^{a-2p}
 =O(n^{(a+1)/p-1})
\]
gives the other two.  All exponents are negative because
$0<\alpha<1$, $p>2$, and $1<a<p-1$; hence
$\sqrt n\,\|\mu^{(n)}-\mu\|_w\to0$ in both cases.

The summability conditions in \eqref{eq:main-countable-summability}
follow from the same parameter ranges.  Since
$\mu_i^{(n)}\leq\mu_i$ for every $i\geq2$, for $m\geq1$,
\[
 \sup_n\sum_{i>m}w_i\mu_i^{(n)}
 \leq\sum_{i>m}w_i\mu_i\longrightarrow0.
\]
This verifies \eqref{eq:main-countable-counts} with $\beta=0$.
The fixed kernel has $\Lambda=0$, is bounded and strictly positive, and
is supercritical because it dominates the constant kernel $\lambda>1$.
Theorem~\ref{thm:countable-concrete} therefore gives the stated limits.
For $\varepsilon>0$, the diagonal choice gives an operator with infinitely
many nonzero diagonal entries plus a rank one operator, and hence has
infinite rank.
\end{proof}

\subsection{Finite types with random parameters}\label{sec:finite-random-proof}

For finitely many types, the deterministic limit is uniform over bounded
scaled parameter perturbations.  We prove this uniformity first, then
condition on the random partition and kernel.

Fix the limiting finite type model of Corollary~\ref{cor:finite-random},
and use $w_i=1$ in this subsection.
For a deterministic partition and kernel $e$, define
\begin{equation}
 \beta_n(e)=\sqrt n\{\mu^{(n)}(e)-\mu\},\qquad
 \Lambda_n(e)=\sqrt n\{\kappa_n(e)-\kappa\}.
\label{eq:finite-environment-perturbations}
\end{equation}
These are the values of $H_n$ and $A_n$ in a fixed environment,
as defined in \eqref{eq:countable-iid-empirical-fluctuation} and
\eqref{eq:finite-parameter-fluctuations}.
Let $\mathfrak E_n(R)$ be the set of partitions and symmetric nonnegative
kernels satisfying the graph model and
$\|\beta_n(e)\|_2+\|\Lambda_n(e)\|_{\mathrm F}\leq R$.
Write $\E_{n,e}$ for expectation over the edges with $e$ fixed, and put
\begin{equation}
\begin{gathered}
 d(h,a)=J^{-1}\{D_\rho h+D_{\mu(1-\rho)}a(\mu\rho)\},\\
 \Sigma=J^{-1}D_{\mu\rho(1-\rho)}J^{-\mathsf T},\\
 V_{n,e}=\frac{L_n-n\mu\rho}{\sqrt n}
             -d\bigl(\beta_n(e),\Lambda_n(e)\bigr).
\end{gathered}
\label{eq:finite-environment-shift}
\end{equation}
The vector $d$ is the parameter shift, while $V_{n,e}$ isolates the edge
fluctuation.

\begin{lemma}[Uniform finite type limit]\label{lem:finite-uniform-law}
For every fixed $R>0$ and $\theta\in\R^K$,
\begin{equation}
 \sup_{e\in\mathfrak E_n(R)}
 \left|\E_{n,e}e^{\mathrm i\theta^{\mathsf T}V_{n,e}}
       -e^{-\theta^{\mathsf T}\Sigma\theta/2}\right|
 \longrightarrow0.
 \label{eq:finite-uniform-centered-law}
\end{equation}
\end{lemma}

\begin{proof}
If \eqref{eq:finite-uniform-centered-law} failed, there would be a
subsequence of environments keeping the error bounded away from zero.
Finite dimensional compactness gives a further subsequence along which
$(\beta_n(e),\Lambda_n(e))\to(\beta,\Lambda)$.
The finite type case of Theorem~\ref{thm:countable-concrete}, with
$w_i=1$, applies along this subsequence: its summability and tail
assumptions are automatic, and the kernels are uniformly bounded.
It gives $V_{n,e}\Rightarrow\Normal_K(0,\Sigma)$, contradicting the
chosen error.
\end{proof}

\begin{proof}[Proof of Corollary~\ref{cor:finite-random}]
Set $V_n=(L_n-n\mu\rho)/\sqrt n-d(H_n,A_n)$.  Conditional on
$\mathcal E_n$, this is $V_{n,e}$ at the realized environment.
For every fixed $\theta\in\R^K$, Lemma~\ref{lem:finite-uniform-law}
and tightness of $(H_n,A_n)$ imply
\[
 \E\left|\E[e^{\mathrm i\theta^{\mathsf T}V_n}\mid\mathcal E_n]
              -e^{-\theta^{\mathsf T}\Sigma\theta/2}\right|\longrightarrow0.
\]
For each $R$, the expectation is bounded by the supremum in
\eqref{eq:finite-uniform-centered-law} plus
$2\Pp(\|H_n\|_2+\|A_n\|_{\mathrm F}>R)$; let $n\to\infty$ and then
$R\to\infty$.
Multiplication by
$\exp\{\mathrm i(u^{\mathsf T}H_n+\operatorname{tr}(B^{\mathsf T}A_n))\}$,
for $u\in\R^K$ and $B\in\operatorname{Sym}_K$, now gives factorization
of the joint characteristic functions.
Since $(H_n,A_n)\Rightarrow(H,A)$, L\'evy's theorem yields
\[
 (H_n,A_n,V_n)
 \Rightarrow\operatorname{Law}(H,A)\boxtimes\Normal_K(0,\Sigma).
\]
The product law makes the last coordinate independent of the parameter
limit.  Applying the
continuous map $(h,a,u)\mapsto d(h,a)+u$ gives a limit whose conditional
law, given $(H,A)$, is $\Normal_K(d(H,A),\Sigma)$.
This is the Gaussian mixture in \eqref{eq:finite-random-limit}.
\end{proof}

\section{Proofs for types on a closed manifold}\label{sec:compact}

Finite spectral projections identify the Gaussian exploration noise;
uniform spectral tail bounds give convergence in $H^{-s}$. Together
with the inverse of $J$, this yields Theorem~\ref{thm:compact} through
\eqref{eq:linearized-exploration}. For independently sampled types,
Corollary~\ref{cor:compact-iid} also requires the joint limit of this
noise and the empirical type fluctuation.

\subsection{Sobolev estimates and invertibility of \texorpdfstring{$J$}{J}}\label{sec:compact-analysis}

The point mass bound makes the exploration noise an $H^{-s}$ valued
random element, and its uniform spectral tail will give tightness.
Multiplication by $C^s$ functions and smoothing by the kernel make the
residual maps well defined; the exponential remainder then gives their
Fr\'echet derivatives.  Finally, comparison with the branching operator
outside the giant gives the inverse of $J$.

Put $a_j=(1+\lambda_j)^{-s}$.  The spectral projections and the
largest discarded point-mass tails are
\begin{equation}
 P_m\eta=\sum_{j=0}^m\eta(e_j)e_j\mathrm{vol}_g,
 \qquad t_m=\sup_{x\in M}\sum_{j>m}a_je_j(x)^2,
 \qquad C_s=\sup_{x\in M}\|\delta_x\|_{-s}.
 \label{eq:compact-spectral-projection}
\end{equation}
For $d=0$, projections beyond the finite basis are the identity.
The next lemma proves $C_s<\infty$ and $t_m\to0$; these are the bounds
needed to control the spectral tails of the exploration noise.

\begin{lemma}[Point masses and weak convergence]\label{lem:compact-analysis}
For the projections and tails in \eqref{eq:compact-spectral-projection},
point masses satisfy
\begin{equation}
 C_s<\infty,\qquad
 \sup_{x\in M}\|(I-P_m)\delta_x\|_{-s}^2=t_m\longrightarrow0.
 \label{eq:compact-delta-tail}
\end{equation}
On positive measures of uniformly bounded mass, weak convergence is
equivalent to convergence in $H^{-s}$.
\end{lemma}

\begin{proof}
For $d=0$, the assertions follow from norm equivalence in a finite
dimensional space.  Suppose that $d\geq1$ and choose
$0<\alpha<\min\{1,s-d/2\}$.  On each connected component, let
$C^{0,\alpha}$ carry the supremum norm plus the $\alpha$-H\"older seminorm,
using its Riemannian distance $\operatorname{dist}_g$.
Sobolev embedding, as in \cite[Chapter~2]{Hebey2000} and
\cite[Chapter~6]{Canzani2013}, gives on each component
\[
 \|f\|_{C^{0,\alpha}}\leq C\|f\|_s,
 \qquad
 \|\delta_x-\delta_y\|_{-s}
 =\sup_{\|f\|_s\leq1}|f(x)-f(y)|
 \leq C\operatorname{dist}_g(x,y)^\alpha.
\]
A compact manifold has finitely many components. Thus
$x\mapsto\delta_x$ is continuous from $M$ to $H^{-s}$ and has compact
image, so $C_s<\infty$.
The orthogonal projections $P_m$ converge strongly to the identity and have
norm at most one.  For a finite $\varepsilon$ net
$\{\delta_{x_1},\ldots,\delta_{x_N}\}$ in that image,
\[
 \sup_x\|(I-P_m)\delta_x\|_{-s}
 \leq\varepsilon+\max_{1\leq\ell\leq N}
                  \|(I-P_m)\delta_{x_\ell}\|_{-s}.
\]
First let $m\to\infty$ and then $\varepsilon\downarrow0$.  This proves the
uniform tail limit in \eqref{eq:compact-delta-tail}, whose spectral
expression is $t_m$.

If $\sigma$ is a positive measure, integration of point masses gives
\[
 \|(I-P_m)\sigma\|_{-s}
 \leq\int_M\|(I-P_m)\delta_x\|_{-s}\,\sigma(\dd x)
 \leq\sigma(M)t_m^{1/2}.
\]
Consequently, if $\sigma_n\Rightarrow\sigma$ weakly and the masses are
bounded, finite spectral projections converge and the displayed inequality
controls their tails.  It follows that
$\|\sigma_n-\sigma\|_{-s}\to0$.  Conversely, convergence in $H^{-s}$ gives
convergence against every smooth function.  Smooth functions are uniformly
dense in $C(M)$; approximation in supremum norm, together with the mass
bound, gives convergence against every continuous function.
\end{proof}

The support of $\mu$ may be a proper compact subset of $M$; the
smooth kernels are still defined on the ambient product $M\times M$.  The Sobolev order and kernel regularity refer to
this ambient manifold.
Among the integer orders allowed here, the smallest is
$s=\lfloor d/2\rfloor+1$.  A smaller $s$ gives a stronger $H^{-s}$
topology and requires fewer kernel derivatives; a larger $s$ gives the
opposite tradeoff.  The present proof uses $C^{2s}$ kernels throughout.
The Sobolev norms are unchanged by rotating an orthonormal basis within
an eigenspace of $-\Delta_g$.  The finite projections used in the proof
may depend on that choice, but the tail estimates hold for every such
basis and the limiting laws do not depend on it.

For $\eta\in H^{-s}$, the factor $e^{-\mathsf K\eta}$ in
\eqref{eq:exploration-residuals} must act boundedly on $H^{-s}$.
The next lemma provides both the smoothing of $\mathsf K\eta$ and the
multiplier bound needed for this product.

\begin{lemma}[Multiplication and kernel smoothing]\label{lem:compact-calculus}
We have
\begin{equation}
 \|af\|_s\leq C\|a\|_{C^s}\|f\|_s,
 \qquad
 \|a\eta\|_{-s}\leq C\|a\|_{C^s}\|\eta\|_{-s}
 \label{eq:compact-multiplier}
\end{equation}
for $a\in C^s(M)$, $f\in H^s$, and $\eta\in H^{-s}$, and
\begin{equation}
 \|\mathsf K_k\eta\|_{C^s}\leq C\|k\|_{C^{2s}}\|\eta\|_{-s}
 \label{eq:compact-smoothing}
\end{equation}
for $k\in C^{2s}(M\times M)$.
\end{lemma}

\begin{proof}
First suppose that $f$ is smooth.  The integer Sobolev norm is equivalent
to the norm of covariant derivatives through order $s$.  Leibniz's rule
gives, for $0\leq\ell\leq s$,
\[
 \|\nabla^\ell(af)\|_{L^2(\mathrm{vol}_g)}
 \leq C_\ell\sum_{r=0}^{\ell}
       \|\nabla^r a\|_\infty\|\nabla^{\ell-r}f\|_{L^2(\mathrm{vol}_g)}.
\]
Summing over $\ell$ proves the first bound in
\eqref{eq:compact-multiplier}; density of smooth functions extends it to
$f\in H^s$.  For the second, use the dual pairing:
\[
 \|a\eta\|_{-s}
 =\sup_{\|f\|_s\leq1}|\eta(af)|
 \leq C\|a\|_{C^s}\|\eta\|_{-s}.
\]
Finally, work in a fixed finite smooth atlas with a subordinate partition
of unity.  For $k\in C^{2s}$, the map $x\mapsto k(x,\cdot)$ is $C^s$
with values in $H^s$: its coordinate derivatives through order $s$ are
continuous as $H^s$ valued maps, since their derivatives in the second variable
through order $s$ are continuous.  Thus, for a coordinate multiindex
$\gamma$ of order at most $s$,
$\partial_x^\gamma(\mathsf K_k\eta)(x)
=\eta(\partial_x^\gamma k(x,\cdot))$.  Duality gives
\[
 \sup_x|\partial_x^\gamma(\mathsf K_k\eta)(x)|
 \leq\|\eta\|_{-s}\sup_x
                  \|\partial_x^\gamma k(x,\cdot)\|_s
 \leq C\|\eta\|_{-s}\|k\|_{C^{2s}}.
\]
The last bound uses at most $s$ derivatives in either variable, hence at
most $2s$ derivatives in total.  Norm equivalence on the fixed atlas
proves \eqref{eq:compact-smoothing}.
\end{proof}

The fixed point equation supplies the smooth multiplier $\rho$ used in
$J$ and in the parameter drift.  Its construction does not require a
density for the type distribution.

\begin{lemma}[Smooth survival probabilities]\label{lem:compact-survival-regularity}
The branching survival probability in \eqref{eq:compact-fixed-point}
has a $C^{2s}(M)$ representative for which that equation holds at every
point of $M$.  It satisfies
\[
 1-\rho\geq e^{-\|\kappa\|_\infty}>0,
 \qquad 0<\rho<1\quad\mu\hbox{-almost everywhere}.
\]
In particular, $J$ in \eqref{eq:compact-J} is bounded on $H^{-s}$.
\end{lemma}

\begin{proof}
Start with any $\mu$-almost everywhere representative of the branching
survival probability and replace it by
\[
 x\longmapsto1-\exp\{-(\mathsf K(\rho\mu))(x)\},\qquad x\in M.
\]
This leaves $\rho\mu$ unchanged, so the fixed point equation
\eqref{eq:compact-fixed-point} then holds everywhere.  Since $\rho\mu$
is a finite measure, differentiation in local coordinates gives
\[
 \partial_x^\gamma(\mathsf K(\rho\mu))(x)
 =\int_M\partial_x^\gamma\kappa(x,y)\rho(y)\,\mu(\dd y),
 \qquad |\gamma|\leq2s.
\]
The derivatives are continuous by uniform continuity of the kernel
derivatives on $M\times M$.  Only the first kernel variable is
differentiated, so no density of $\mu$ is required.  Thus
$\rho\in C^{2s}(M)$.  Nonnegativity and boundedness of the kernel give
the lower bound on $1-\rho$; irreducibility and supercriticality give
$0<\rho<1$ $\mu$-almost everywhere.  Finally,
\eqref{eq:compact-multiplier}--\eqref{eq:compact-smoothing} show that $J$
is bounded.
\end{proof}

The point mass bounds also ensure that the proposed noise covariances
are trace class operators on $H^{-s}$.

\begin{lemma}[Covariance operators for Borel type laws]\label{lem:compact-covariance-operators}
For every bounded nonnegative Borel function $a$, the operator
$\mathcal V_\mu(a)$ in \eqref{eq:measure-covariance-notation} is positive
and trace class, with
\begin{equation}
 \operatorname{tr}\mathcal V_\mu(a)
 =\int_M a(x)\|\delta_x\|_{-s}^2\,\mu(\dd x)
 \leq C_s^2\int_M a\,\dd\mu.
 \label{eq:compact-noise-operator}
\end{equation}
The Gaussian law $\Normal(0,\mathcal V_\mu(a))$ has covariance under
evaluation at $f,g\in H^s$ given by
\begin{equation}
 (f,g)\longmapsto\int_M afg\,\dd\mu.
 \label{eq:compact-noise}
\end{equation}
\end{lemma}

\begin{proof}
The continuity of $x\mapsto\delta_x$ in
Lemma~\ref{lem:compact-analysis} makes
$x\mapsto\delta_x\otimes\delta_x$ continuous in trace norm.  Multiplication
by $a$ gives a strongly measurable integrand, and the rank one identity
\[
 \|\delta_x\otimes\delta_x\|_{\mathrm{tr}}
 =\|\delta_x\|_{-s}^2\leq C_s^2
\]
proves Bochner integrability in the trace class operators on $H^{-s}$.
Positivity and \eqref{eq:compact-noise-operator} follow by integration.
Applying the continuous evaluation functionals at $f,g$ gives
\eqref{eq:compact-noise}.
\end{proof}

This construction applies to every Borel probability measure $\mu$,
including measures supported on lower dimensional submanifolds.

The derivative comparison in \eqref{eq:compact-derivative} and the
parameter expansion in \eqref{eq:compact-drift} require a quadratic
remainder for the exponential. The product bound in $C^s$ gives, on a
fixed ball of radius $R$,
\begin{equation}
 \|e^a-e^b\|_{C^s}\leq C_R\|a-b\|_{C^s},\qquad
 \|e^{a+b}-e^a-e^ab\|_{C^s}\leq C_R\|b\|_{C^s}^2
 \label{eq:compact-exponential}
\end{equation}
when $\|a\|_{C^s}\leq R$ and $\|b\|_{C^s}\leq R$ in the second bound.
Indeed, the identities
\[
 e^{a+b}-e^a=\int_0^1e^{a+tb}b\,\dd t,\qquad
 e^{a+b}-e^a-e^ab=\int_0^1(1-t)e^{a+tb}b^2\,\dd t
\]
hold in $C^s$ and give these estimates by Leibniz's rule.

The factorization used in Lemma~\ref{lem:countable-inverse} has a
continuous type version: the types outside the giant again determine a
symmetric branching operator. Its spectral gap will give the inverse of
$J$ on $H^{-s}$.

\begin{lemma}[Invertibility of the limiting derivative]
\label{lem:compact-inverse}
Suppose that $\kappa\in C^{2s}(M\times M)$ is symmetric, nonnegative,
irreducible with respect to $\mu$, and
$\|T_\kappa\|_{L^2(\mu)\to L^2(\mu)}>1$.
Then the operator $J$ in \eqref{eq:compact-J} is boundedly invertible on
$H^{-s}$.
\end{lemma}

\begin{proof}
\noindent Step 1: The symmetric branching operator.
Put $\nu=(1-\rho)\mu$.  The measure $\nu$ describes the limiting types outside
the giant. Define $U:L^2(\nu)\to H^{-s}$ by $Uh=h\nu$, where $h\nu$ is
the signed measure with density $h$ relative to $\nu$.
The point mass bound \eqref{eq:compact-delta-tail} gives
\[
 \|Uh\|_{-s}\leq C_s\int_M|h|\,\dd\nu
 \leq C_s\nu(M)^{1/2}\|h\|_{L^2(\nu)}.
\]
Thus $\|U\|\leq C_s\nu(M)^{1/2}$, and
\eqref{eq:compact-smoothing} gives
$\|\mathsf K\|_{H^{-s}\to L^2(\nu)}
\leq C\|\kappa\|_{C^{2s}}\nu(M)^{1/2}$.  The operator $A=\mathsf KU$ on $L^2(\nu)$
is the integral operator with kernel $\kappa$.  It is selfadjoint and
Hilbert--Schmidt because
\[
 \iint\kappa(x,y)^2\,\nu(\dd x)\nu(\dd y)
 \leq\|\kappa\|_\infty^2\nu(M)^2<\infty.
\]
Its operator norm $r=\|A\|_{L^2(\nu)\to L^2(\nu)}$ is a positive eigenvalue with a nonnegative eigenfunction
$g\ne0$.  To see this, compact selfadjointness gives an eigenfunction at an
eigenvalue of absolute value $r$.  Since the kernel is nonnegative,
$|\langle f,Af\rangle_{L^2(\nu)}|
\leq\langle |f|,A|f|\rangle_{L^2(\nu)}$.
Thus the largest Rayleigh quotient is $r$ and is attained by a nonnegative
function.  Here $r>0$ because supercriticality makes the kernel nonzero
and $1-\rho\geq e^{-\|\kappa\|_\infty}$.

\par\smallskip\noindent Step 2: Subcriticality.
Compare the action of this operator on the
positive function $t=\rho/(1-\rho)$.  This function is bounded and strictly positive $\nu$ almost
everywhere.  The fixed point equation \eqref{eq:compact-fixed-point} gives
\[
 At=\mathsf K(\rho\mu)=-\log(1-\rho)<\frac{\rho}{1-\rho}=t
 \quad\nu\hbox{ almost everywhere}.
\]
Here $-(1-\rho)\log(1-\rho)<\rho$ whenever $0<\rho<1$. Pairing with $g$ gives
\[
 r\langle g,t\rangle_{L^2(\nu)}
 =\langle g,At\rangle_{L^2(\nu)}
 <\langle g,t\rangle_{L^2(\nu)}.
\]
The last inner product is strictly positive, so $r<1$.  Hence
$(I-A)^{-1}=\sum_{k\geq0}A^k$ is a bounded operator on $L^2(\nu)$.
\par\smallskip\noindent Step 3: The inverse on the Sobolev space.
To obtain the inverse on $H^{-s}$, note that $J=I-U\mathsf K$.
The identity for the inverses of $I-U\mathsf K$ and $I-\mathsf KU$
suggests the candidate
\[
 B:=I+U(I-\mathsf KU)^{-1}\mathsf K.
\]
Direct multiplication gives $JB=BJ=I$.  The Neumann series for $I-A$
therefore gives the quantitative bound
\begin{equation}
 \|J^{-1}\|_{\mathrm{op}}
 \leq1+\frac{\|U\|_{L^2(\nu)\to H^{-s}}
                   \|\mathsf K\|_{H^{-s}\to L^2(\nu)}}{1-r}
 \leq1+\frac{C_sC\|\kappa\|_{C^{2s}}\nu(M)}{1-r}.
 \label{eq:compact-inverse-bound}
\end{equation}
In particular, $J^{-1}=B$ is bounded on $H^{-s}$.
\end{proof}

For finite signed measures $\eta$, Fubini's theorem verifies
\eqref{eq:compact-transpose} directly from \eqref{eq:compact-J}, since
$(J\eta)(f)=\eta(J^*f)$.
The multiplier and smoothing estimates make both sides continuous in
$\eta\in H^{-s}$, so density of smooth densities extends the identity to
that space.  Bounded invertibility of $J$ then gives bounded invertibility
of $J^*$ and the identity
\[
 (J^{-1}\eta)(f)=\eta(Rf),\qquad \eta\in H^{-s},\quad f\in H^s.
\]

\subsection{The limit for deterministic types}

We first approximate the giant by the seeded exploration, then prove the
martingale limit and analyze \eqref{eq:exploration-equation}.  Throughout this
subsection, $M,s,\mu,\rho$ are as in Subsection~\ref{sec:compact-model},
and the limiting kernel $\kappa\in C^{2s}(M\times M)$ is symmetric,
nonnegative, irreducible with respect to $\mu$, and supercritical as in
\eqref{eq:compact-model}.  The types and kernels $\kappa_n$ are deterministic,
and each $\kappa_n\in C^{2s}$ is symmetric and nonnegative.  For the
preparatory lemmas we assume
\begin{equation}
 \mu^{(n)}\Rightarrow\mu\quad\hbox{weakly},\qquad
 \|\kappa_n-\kappa\|_{C^{2s}}\longrightarrow0.
 \label{eq:compact-ordinary-convergence}
\end{equation}
The first order assumptions in \eqref{eq:compact-expansions} are needed
only for the expansion at $\rho\mu$ in
Lemma~\ref{lem:compact-terminal-maps} and in the proof of the theorem.

\begin{lemma}[Approximation by the seeded exploration]
\label{lem:compact-seeded-approximation}
The seeds meet $C_1^{(n)}$ with probability tending to one, and
\begin{equation}
 \|\mathcal R_{n,\tau_n}-\mathcal L_n\|_{-s}/\sqrt n=\op(1),\qquad
 \mathcal R_{n,\tau_n}/n\xrightarrow{\mathrm p}\rho\mu
       \quad\hbox{in }H^{-s}.
 \label{eq:compact-exploration}
\end{equation}
\end{lemma}

\begin{proof}
We apply the component results of Bollob\'as, Janson, and Riordan
\cite[Theorems~9.10 and~3.12(ii)]{BollobasJansonRiordan2007}.
The kernels are continuous and converge uniformly to $\kappa$, and
$\sup_n\|\kappa_n\|_\infty<\infty$.  Writing $e(G_n)$ for the number of
edges, the exponential edge probabilities give
\[
 \frac1n\E e(G_n)
 =\frac12\iint\kappa_n(x,y)\,\mu^{(n)}(\dd x)\mu^{(n)}(\dd y)+O(n^{-1})
 \longrightarrow\frac12\iint\kappa\,\dd\mu\,\dd\mu.
\]
The error follows from the omission of diagonal vertex pairs and the
uniform expansion $n(1-e^{-\kappa_n/n})=\kappa_n+O(n^{-1})$.
For the integral limit, use uniform convergence of the kernels and
$\mu^{(n)}\boxtimes\mu^{(n)}\Rightarrow\mu\boxtimes\mu$ on the compact product
space.  These observations verify graphicality; the exponential
parametrization is covered by
\cite[Remark~2.4]{BollobasJansonRiordan2007}.  The cited theorems yield
\begin{equation}
 \mathcal L_n/n\xrightarrow{\mathrm p}\rho\mu
       \quad\hbox{weakly},\qquad |C_2^{(n)}|=\Op(\log n).
 \label{eq:compact-component-first-order}
\end{equation}
The measures in \eqref{eq:compact-component-first-order} have masses at
most one, so the first convergence holds in
$H^{-s}$ by Lemma~\ref{lem:compact-analysis}.

By \eqref{eq:seed-hit}, with $a=\frac12\int_M\rho\,\dd\mu$, and the size law of large
numbers in \eqref{eq:compact-component-first-order}, the seeds meet the
giant with probability tending to one.  On this event, at most
$\ell_n|C_2^{(n)}|$ additional vertices are explored.  Hence
\eqref{eq:compact-delta-tail} gives
\[
 \|\mathcal R_{n,\tau_n}-\mathcal L_n\|_{-s}
 \leq C_s\ell_n|C_2^{(n)}|=\Op((\log n)^2)
\]
outside an event of probability tending to zero.  This proves the first
assertion in \eqref{eq:compact-exploration}; combining it with
\eqref{eq:compact-component-first-order} and
Lemma~\ref{lem:compact-analysis} proves the second.
\end{proof}

The law of large numbers already determines the survival factors.
By \eqref{eq:compact-smoothing} and \eqref{eq:compact-exploration},
\[
 \|\mathsf K_{\kappa_n}(\mathcal R_{n,\tau_n}/n)-\mathsf K(\rho\mu)\|_{C^s}
 \leq C\|\mathcal R_{n,\tau_n}/n-\rho\mu\|_{-s}
       +C\|\kappa_n-\kappa\|_{C^{2s}}\|\rho\mu\|_{-s}
 \xrightarrow{\mathrm p}0.
\]
Since \eqref{eq:exploration-Q} gives
$Q_{n,\tau_n}=e^{-\mathsf K_{\kappa_n}(\mathcal R_{n,\tau_n}/n)}$
and $(1-\rho)=e^{-\mathsf K(\rho\mu)}$, the exponential estimate
\eqref{eq:compact-exponential} implies
\begin{equation}
 Q_{n,\tau_n}\xrightarrow{\mathrm p}(1-\rho),\qquad
 Q_{n,\tau_n}^{-1}\xrightarrow{\mathrm p}(1-\rho)^{-1}
 \quad\hbox{in }C^s.
 \label{eq:compact-survival-convergence}
\end{equation}
We next study the normalized martingale increment $N_n$ from
\eqref{eq:exploration-martingale-noise}, now viewed as an element of
$H^{-s}$.

\begin{lemma}[The martingale limit in the negative Sobolev space]
\label{lem:compact-martingale-limit}
For $N_n$ in \eqref{eq:exploration-martingale-noise}, we have
\begin{equation}
 \sup_{t\leq\tau_n}\|M_{n,t}-\mu^{(n)}\|_{-s}=\Op(n^{-1/2}),\qquad
 \sup_n\E\|(I-P_m)N_n\|_{-s}^2\leq Ct_m.
 \label{eq:compact-martingale-bounds}
\end{equation}
Moreover,
\begin{equation}
 N_n\Rightarrow\Normal\!\left(0,
       \mathcal V_\mu\!\left(\frac{\rho}{1-\rho}\right)\right)
 \quad\hbox{in }H^{-s}.
 \label{eq:compact-martingale-covariance}
\end{equation}
The constants in the moment estimates are uniform over the type
configurations.
\end{lemma}

\begin{proof}
\noindent Step 1: Second moments and spectral tails.
Bounds for the whole martingale and its discarded spectral coordinates
will allow us to pass from finite spectral Gaussian limits to an $H^{-s}$ limit.
Choose $\kappa_*<\infty$ with $\sup_n\|\kappa_n\|_\infty\leq\kappa_*$.
Equation~\eqref{eq:exploration-Q} gives
$e^{-\kappa_*}\leq Q_{n,t}\leq1$.
In particular, \eqref{eq:compact-delta-tail} gives
$\|W_{n,t}\|_{-s}\leq e^{\kappa_*}C_sn$.
The increment \eqref{eq:measure-increment} is a finite sum of conditionally
independent centered point masses, so $W_{n,\cdot}$ is an adapted
$H^{-s}$ valued square integrable martingale.  For $t<\tau_n$, its
conditional second moment is
\[
 \E(\|\Delta W_{n,t}\|_{-s}^2\mid\mathcal F_{n,t})
 =\sum_{v\in U_{n,t}}
 \frac{p_{n,t}(x_v^{(n)})}
      {Q_{n,t}(x_v^{(n)})^2(1-p_{n,t}(x_v^{(n)}))}
      \|\delta_{x_v^{(n)}}\|_{-s}^2\leq C,
\]
since $p_{n,t}\leq C/n$ and $|U_{n,t}|\leq n$.
Using \eqref{eq:compact-delta-tail}, replacing each point mass by
$(I-P_m)\delta_{x_v^{(n)}}$ gives the bound
$Ct_m$ for the projected conditional second moment.
The martingale isometry and Doob's inequality therefore yield
\[
 \E\sup_{t\leq\tau_n}\|W_{n,t}-W_{n,0}\|_{-s}^2\leq Cn,\qquad
 \E\|(I-P_m)(W_{n,\tau_n}-W_{n,0})\|_{-s}^2\leq Cnt_m.
\]
Because $W_{n,0}=n\mu^{(n)}-S_n$ and
$\|S_n\|_{-s}\leq C_s\ell_n$, this proves
\eqref{eq:compact-martingale-bounds}; the seed contribution to its first
estimate is $O(\ell_n/n)=o(n^{-1/2})$.
The bounds in \eqref{eq:compact-martingale-bounds} depend only on
\eqref{eq:compact-delta-tail} and the uniform kernel bound, so they are
uniform over the type configurations.

\par\smallskip\noindent Step 2: The limiting bracket.
To identify the covariance of each finite projection, take $f,g$ to be
finite linear combinations of the
eigenfunctions $e_j$.
The exact bracket identity \eqref{eq:common-bracket} gives
\begin{align}
 \frac1n\langle W_{n,\cdot}(f),W_{n,\cdot}(g)\rangle_{\tau_n}
 &=\sum_{t<\tau_n}M_{n,t}
       \bigl(fg(Q_{n,t+1}^{-1}-Q_{n,t}^{-1})\bigr)\nonumber\\
 &=\mu^{(n)}\bigl(fg(Q_{n,\tau_n}^{-1}-1)\bigr)+\op(1).
 \label{eq:compact-bracket}
\end{align}
To justify the replacement, note that the positive measure
$\mathcal R_{n,t}/n$ has mass at most one.  Direct integration of kernel
derivatives in \eqref{eq:exploration-Q} gives
$\|Q_{n,t}^{\pm1}\|_{C^s}\leq C$.  For $t<\tau_n$,
\[
 Q_{n,t+1}^{-1}-Q_{n,t}^{-1}
 =Q_{n,t}^{-1}
      \{\exp(\kappa_n(\cdot,x_{\omega_{n,t}}^{(n)})/n)-1\},
 \qquad
 \|Q_{n,t+1}^{-1}-Q_{n,t}^{-1}\|_{C^s}\leq C/n
\]
by \eqref{eq:compact-exponential}.  The multiplier estimate
\eqref{eq:compact-multiplier} and \eqref{eq:compact-martingale-bounds}
bound the replacement error in \eqref{eq:compact-bracket} by
\[
 \sup_{t\leq\tau_n}\|M_{n,t}-\mu^{(n)}\|_{-s}
 \sum_{t<\tau_n}\|fg(Q_{n,t+1}^{-1}-Q_{n,t}^{-1})\|_s
 \leq C_{f,g}\sup_{t\leq\tau_n}\|M_{n,t}-\mu^{(n)}\|_{-s}=\op(1).
\]
By \eqref{eq:compact-survival-convergence}, the functions
$Q_{n,\tau_n}^{-1}$ converge to $(1-\rho)^{-1}$ in $C^s$ in probability.
By Lemma~\ref{lem:compact-analysis}, weak convergence of
$\mu^{(n)}$ implies convergence in $H^{-s}$.  The multiplier estimate
\eqref{eq:compact-multiplier} now
shows that the right side of \eqref{eq:compact-bracket} tends to
$\int fg\frac{\rho}{1-\rho}\,\dd\mu$.

\par\smallskip\noindent Step 3: Gaussian convergence in the Sobolev space.
The conditional Lindeberg bound in Lemma~\ref{lem:exploration-lindeberg}, applied to
finite linear combinations of eigenfunctions, and the scalar martingale
central limit theorem \cite[Corollary~3.1]{HallHeyde1980} give the
Gaussian limit of every finite spectral projection of $N_n$.
The candidate covariance operator is positive and trace class by
\eqref{eq:compact-noise-operator}, with trace at most
$(e^{\kappa_*}-1)C_s^2$.  Its discarded Gaussian second moment is
\[
 \operatorname{tr}\left\{(I-P_m)
       \mathcal V_\mu\!\left(\frac{\rho}{1-\rho}\right)(I-P_m)\right\}
 =\int_M\frac{\rho}{1-\rho}\sum_{j>m}a_je_j(x)^2\,\mu(\dd x)
 \leq Ct_m.
\]
Together with the uniform tail bound
\eqref{eq:compact-martingale-bounds}, the finite spectral limits satisfy
Lemma~\ref{lem:hilbert-projections}.  Hence
\begin{equation}
 N_n\Rightarrow\Normal\!\left(0,
       \mathcal V_\mu\!\left(\frac{\rho}{1-\rho}\right)\right)
 \quad\hbox{in }H^{-s}.
 \label{eq:compact-normalized-noise}
\end{equation}
The covariance identity extends from finite linear combinations of the
eigenfunctions to all $f,g\in H^s$ by density and the bound
$H^s\subset C(M)$.
\end{proof}

The smoothing and multiplier estimates
\eqref{eq:compact-smoothing} and~\eqref{eq:compact-multiplier} extend the
residual maps \eqref{eq:exploration-residuals} continuously to $H^{-s}$.
The next lemma identifies the Gaussian limit of the noise in the common
exploration equation \eqref{eq:exploration-equation}.

\begin{lemma}[Noise in the exploration equation]
\label{lem:compact-terminal-noise}
The exploration noise satisfies
\begin{equation}
 Z_n\Rightarrow\Normal\!\left(0,\mathcal V_\mu(\rho(1-\rho))\right)
 \quad\hbox{in }H^{-s}.
 \label{eq:compact-terminal-noise}
\end{equation}

\end{lemma}

\begin{proof}
Lemma~\ref{lem:compact-martingale-limit} gives tightness of $N_n$,
and \eqref{eq:compact-survival-convergence} gives
$Q_{n,\tau_n}\to (1-\rho)$ in $C^s$.  The multiplier bound
\eqref{eq:compact-multiplier}, applied to
\eqref{eq:exploration-noise}, gives
\begin{equation}
 \|Z_n-(1-\rho)N_n\|_{-s}
 \leq C\|Q_{n,\tau_n}-(1-\rho)\|_{C^s}\|N_n\|_{-s}
      +C\|Q_{n,\tau_n}\|_{C^s}C_s\ell_n/\sqrt n
 \xrightarrow{\mathrm p}0.
 \label{eq:compact-noise-multiplier-error}
\end{equation}
The continuous map $\eta\mapsto(1-\rho)\eta$ transforms the
Gaussian law in \eqref{eq:compact-normalized-noise} into a centered
Gaussian law with evaluation covariance
\[
 (f,g)\longmapsto\int fg(1-\rho)^2\frac{\rho}{1-\rho}\,\dd\mu
 =\int fg(1-\rho)\rho\,\dd\mu.
\]
This is the covariance of $\Normal(0,\mathcal V_\mu(\rho(1-\rho)))$.
Slutsky's theorem and \eqref{eq:compact-noise-multiplier-error} prove the
asserted convergence.
\end{proof}

The derivative convergence and first order value of $h_n(\rho\mu)$ are the
two remaining hypotheses of Lemma~\ref{lem:linearization}.

\begin{lemma}[Differentiability and first order expansion]
\label{lem:compact-terminal-maps}
The maps $h_n,h:H^{-s}\to H^{-s}$ are continuously Fr\'echet
differentiable, $h(\rho\mu)=0$, and $Dh(\rho\mu)=J$.
For every fixed $R>0$,
\begin{equation}
 \sup_{\|\eta-\rho\mu\|_{-s}\leq R}
      \|Dh_n(\eta)-Dh(\eta)\|_{\mathrm{op}}
 \leq C_R\{\|\mu^{(n)}-\mu\|_{-s}
                    +\|\kappa_n-\kappa\|_{C^{2s}}\}\longrightarrow0.
 \label{eq:compact-derivative}
\end{equation}
If, in addition, the first order expansions
\eqref{eq:compact-expansions} hold, then
\begin{equation}
 \sqrt n\,h_n(\rho\mu)\longrightarrow-b,\qquad
 b:=\rho\beta+(1-\rho)(\mathsf K_\Lambda(\rho\mu))\mu
 \quad\hbox{in }H^{-s}.
 \label{eq:compact-drift}
\end{equation}
\end{lemma}

\begin{proof}
Write $F_n(\eta)=e^{-\mathsf K_{\kappa_n}\eta}$ and
$F(\eta)=e^{-\mathsf K\eta}$. The fixed point equation \eqref{eq:compact-fixed-point} and
\eqref{eq:exploration-residuals} give $h(\rho\mu)=0$.
On a fixed ball $\|\eta-\rho\mu\|_{-s}\leq R$,
\eqref{eq:compact-smoothing} bounds the
exponential arguments in $C^s$.  The Taylor estimate
\eqref{eq:compact-exponential} and \eqref{eq:compact-multiplier}
give, for
$\|v\|_{-s}\leq1$,
\[
 \|h_n(\eta+v)-h_n(\eta)
          -v+e^{-\mathsf K_{\kappa_n}\eta}
                    (\mathsf K_{\kappa_n}v)\mu^{(n)}\|_{-s}
 \leq C_R\|v\|_{-s}^2.
\]
Here $\|\mu^{(n)}\|_{-s}\leq C_s$.  Thus $h_n$ is Fr\'echet differentiable,
with
\[
 Dh_n(\eta)v
 =v-e^{-\mathsf K_{\kappa_n}\eta}(\mathsf K_{\kappa_n}v)\mu^{(n)}.
\]
The same calculation applies to $h$ with $\mu,\kappa$ in place of
$\mu^{(n)},\kappa_n$; at $\rho\mu$ it gives $Dh(\rho\mu)=J$.
The derivative formulas are locally Lipschitz in $\eta$ by
\eqref{eq:compact-exponential} and \eqref{eq:compact-smoothing}, hence
continuous.

For the derivative comparison, if $\|v\|_{-s}\leq1$, the difference
$Dh_n(\eta)v-Dh(\eta)v$ is the negative sum of
\[
 \{F_n(\eta)-F(\eta)\}(\mathsf K_{\kappa_n}v)\mu^{(n)},
 \quad F(\eta)(\mathsf K_{\kappa_n}-\mathsf K)v\,\mu^{(n)},
 \quad F(\eta)(\mathsf Kv)(\mu^{(n)}-\mu).
\]
The first two terms are bounded in $H^{-s}$ by
$C_R\|\kappa_n-\kappa\|_{C^{2s}}\|v\|_{-s}$, and the third by
$C_R\|\mu^{(n)}-\mu\|_{-s}\|v\|_{-s}$.
Taking the supremum over $\|v\|_{-s}\leq1$ and
$\|\eta-\rho\mu\|_{-s}\leq R$ proves \eqref{eq:compact-derivative}.

It remains to determine $\sqrt n\,h_n(\rho\mu)$.  Under
\eqref{eq:compact-expansions}, put $d_n=\mathsf K_{\kappa_n-\kappa}(\rho\mu)$.
This is the change in $\mathsf K_\kappa(\rho\mu)$ caused by the kernel
perturbation.  The kernel expansion gives
\[
 \|d_n\|_{C^s}=O(n^{-1/2}),\qquad
 \sqrt n\,d_n\longrightarrow\mathsf K_\Lambda(\rho\mu)\quad\hbox{in }C^s.
\]
Since $e^{-\mathsf K(\rho\mu)}=(1-\rho)$, the exact expansion is
\[
 h_n(\rho\mu)=-\rho(\mu^{(n)}-\mu)-(1-\rho)d_n\mu+\epsilon_n,
 \quad
 \epsilon_n=(1-\rho)(e^{-d_n}-1+d_n)\mu^{(n)}-(1-\rho)d_n(\mu^{(n)}-\mu).
\]
The estimates \eqref{eq:compact-exponential} and
\eqref{eq:compact-multiplier} give
\[
 \|\epsilon_n\|_{-s}
 \leq C\{\|d_n\|_{C^s}^2+
             \|d_n\|_{C^s}\|\mu^{(n)}-\mu\|_{-s}\}=O(n^{-1}).
\]
Multiplying the expansion by $\sqrt n$ proves \eqref{eq:compact-drift}.
\end{proof}

\begin{proof}[Proof of Theorem~\ref{thm:compact}]
The exploration equation determines the fluctuation once its random error
and its response to a change in the explored mass are known.  We first
locate the explored mass near $\rho\mu$, then linearize that equation and
transfer the result to the giant.

\noindent Step 1: Approximate the giant by the exploration.
Assumption~\assref{C} gives $\|\mu^{(n)}-\mu\|_{-s}\to0$ and
$\|\kappa_n-\kappa\|_{C^{2s}}\to0$.  By
Lemma~\ref{lem:compact-analysis}, the first convergence also holds weakly
as convergence of positive measures; the second implies uniform kernel
convergence.  Thus \eqref{eq:compact-ordinary-convergence} holds, and
Lemma~\ref{lem:compact-seeded-approximation} gives
\[
 \widehat\mu_n\xrightarrow{\mathrm p}\rho\mu,
 \qquad
 \frac{\mathcal R_{n,\tau_n}-\mathcal L_n}{\sqrt n}
 \xrightarrow{\mathrm p}0
 \quad\hbox{in }H^{-s}.
\]
The first convergence selects the nonzero fixed point at which to expand;
the second makes the extra components reached by the seeds negligible on
the fluctuation scale.

\par\smallskip\noindent Step 2: Separate edge noise from parameter drift.
The exact exploration equation \eqref{eq:exploration-equation} and
Lemma~\ref{lem:compact-terminal-noise} give
\[
 \sqrt n\,h_n(\widehat\mu_n)=-Z_n,
 \qquad
 Z_n\Rightarrow\Normal\!\left(0,\mathcal V_\mu(\rho(1-\rho))\right)
 \quad\hbox{in }H^{-s}.
\]
At the limiting mass, Lemma~\ref{lem:compact-terminal-maps} gives
\[
 \sqrt n\,h_n(\rho\mu)\longrightarrow-b,
 \qquad
 b=\rho\beta+(1-\rho)(\mathsf K_\Lambda(\rho\mu))\mu.
\]
The two terms in $b$ are the first order changes caused by the type and
kernel perturbations.  These formulas verify
$\sqrt n\,h_n(\widehat\mu_n)=\Op(1)$ and
$\sqrt n\,h_n(\rho\mu)=O(1)$, the two residual bounds required for local
linearization.

\par\smallskip\noindent Step 3: Linearize and solve for the explored mass.
Lemma~\ref{lem:compact-terminal-maps} gives continuous Fr\'echet
differentiability, $Dh(\rho\mu)=J$, and local uniform derivative
convergence in \eqref{eq:compact-derivative}.
Lemma~\ref{lem:compact-inverse} gives the bounded inverse of $J$.
Therefore Lemma~\ref{lem:linearization} applies with
$g_n=h_n$, $g=h$, $x_\star=\rho\mu$, and $X_n=\widehat\mu_n$.
It first yields $\sqrt n(\widehat\mu_n-\rho\mu)=\Op(1)$ and then
\[
 \sqrt n\{h_n(\widehat\mu_n)-h_n(\rho\mu)\}
 =J\sqrt n(\widehat\mu_n-\rho\mu)+\op(1).
\]
Substituting the residual values from Step~2 and applying $J^{-1}$ gives
\[
 \sqrt n(\widehat\mu_n-\rho\mu)
 =J^{-1}(b-Z_n)+\op(1)
 \quad\hbox{in }H^{-s},
\]
which is \eqref{eq:linearized-exploration} for this model.

\par\smallskip\noindent Step 4: Return to the giant and identify its Gaussian law.
Since $n\widehat\mu_n=\mathcal R_{n,\tau_n}$, the approximation in Step~1
gives
\[
 \frac{\mathcal L_n-n\rho\mu}{\sqrt n}
 =\sqrt n(\widehat\mu_n-\rho\mu)
   +\frac{\mathcal L_n-\mathcal R_{n,\tau_n}}{\sqrt n}
 =J^{-1}(b-Z_n)+\op(1).
\]
The Gaussian limit in Step~2 and Slutsky's theorem give mean $J^{-1}b$
and covariance operator
\begin{equation}
 \mathcal C=J^{-1}\mathcal V_\mu(\rho(1-\rho))
                  (J^{-1})^{\dagger_{-s}}.
 \label{eq:compact-limit-covariance}
\end{equation}
This operator is positive and trace class: by
\eqref{eq:compact-noise-operator} and boundedness of $J^{-1}$,
\[
 \operatorname{tr}\mathcal C
 \leq\|J^{-1}\|_{\mathrm{op}}^2
       \operatorname{tr}\mathcal V_\mu(\rho(1-\rho))<\infty.
\]
Thus the limit is $\Normal(J^{-1}b,\mathcal C)$.
Substitution of $b$ proves \eqref{eq:compact-CLT}; the dual identity
following \eqref{eq:compact-transpose} also gives the scalar formulas
\eqref{eq:compact-scalar-mean}--\eqref{eq:compact-scalar-covariance}.
\end{proof}

\subsection{Independently sampled types}

Random types create a second source of Gaussian fluctuations.  We first
separate this empirical type noise from the edge noise, and then insert
both into the equation for the explored type measure.

Throughout this subsection, the types $X_1,X_2,\ldots$ are independent
with law $\mu$, and the kernel is the fixed $\kappa$ of
Corollary~\ref{cor:compact-iid}.  We use $H_n$ from \eqref{eq:countable-iid-empirical-fluctuation}.
This is the fluctuation of the empirical type distribution
before any edges are drawn. The exploration noise $Z_n$ in
\eqref{eq:exploration-noise} represents the additional randomness from the
edges and the negligible seed correction.

\begin{lemma}[Joint limit of type and edge noise]\label{lem:compact-joint-noise}
In $H^{-s}\times H^{-s}$,
\[
 (H_n,Z_n)\Rightarrow\Normal\!\left(0,
 \begin{pmatrix}
  \mathcal V_\mu(1)-\mu\otimes\mu&0\\
  0&\mathcal V_\mu(\rho(1-\rho))
 \end{pmatrix}\right).
\]
In particular, the two coordinates are asymptotically independent.
\end{lemma}

\begin{proof}
Every finite spectral projection satisfies the
multivariate iid central limit theorem.  Independence and centering give
\begin{equation}
 \E\|(I-P_m)H_n\|_{-s}^2
 =\E\|(I-P_m)(\delta_{X_1}-\mu)\|_{-s}^2
 \leq\E\|(I-P_m)\delta_{X_1}\|_{-s}^2\leq t_m.
 \label{eq:compact-iid-type-tail}
\end{equation}
The candidate empirical covariance is
$\mathcal V_\mu(1)-\mu\otimes\mu
 =\E[(\delta_{X_1}-\mu)\otimes(\delta_{X_1}-\mu)]$.
It is positive and trace class, with trace
$\E\|\delta_{X_1}-\mu\|_{-s}^2\leq C_s^2$ and evaluation covariance
$\mu(fg)-\mu(f)\mu(g)$.
Lemma~\ref{lem:hilbert-projections}, applied to the finite spectral
limits and \eqref{eq:compact-iid-type-tail}, therefore gives
\[
 H_n\Rightarrow\Normal(0,\mathcal V_\mu(1)-\mu\otimes\mu)
 \quad\hbox{in }H^{-s}.
\]

With the kernel fixed, the environment $\mathcal E_n$ of
Subsection~\ref{sec:exploration} is generated by $X_1,\ldots,X_n$.
We realize conditional expectation by integrating the edge and seed
variables while holding the types fixed.  For any deterministic type sequence whose empirical measures converge
weakly to $\mu$, Lemmas~\ref{lem:compact-seeded-approximation}--\ref{lem:compact-terminal-noise}
give convergence of $Z_n$ to
$\Normal(0,\mathcal V_\mu(\rho(1-\rho)))$ with the fixed kernel $\kappa$.
The strong law applied to a countable uniformly dense subset of $C(M)$
gives weak empirical convergence almost surely. For every fixed finite
linear combination $g$ of the eigenfunctions, the conditional
characteristic function satisfies
\[
 \phi_n(g):=\E[e^{\mathrm i Z_n(g)}\mid\mathcal E_n]
 \longrightarrow
 \exp\{-\tfrac12\int g^2(1-\rho)\rho\,\dd\mu\}
 \quad\hbox{almost surely}.
\]
For another such linear combination $f$, conditioning gives
\[
 \E e^{\mathrm iH_n(f)+\mathrm iZ_n(g)}
 =\E[e^{\mathrm iH_n(f)}\phi_n(g)].
\]
Since $|\phi_n(g)|\leq1$, bounded convergence and the finite dimensional
limit of $H_n$ show that this tends to
\[
 \exp\left\{-\frac12\left[\mu(f^2)-\mu(f)^2
                    +\mu(g^2\rho(1-\rho))\right]\right\}.
\]
Finite linear combinations of test functions therefore prove joint finite
dimensional Gaussian convergence with zero cross covariance.

The spectral tail bound \eqref{eq:compact-martingale-bounds} is uniform
in the type configuration, so it also gives annealed tightness of the
normalized increment $N_n$ in \eqref{eq:exploration-martingale-noise}.
For almost every type sequence, \eqref{eq:compact-survival-convergence}
gives $Q_{n,\tau_n}\to (1-\rho)$ in conditional probability in $C^s$.  Integrating
these conditional probabilities gives convergence in probability under
the joint law of types, edges, and seeds.  The estimate
\eqref{eq:compact-noise-multiplier-error} therefore applies.
Thus $Z_n$ is tight, as is $H_n$.  Applying
Lemma~\ref{lem:hilbert-projections} on $H^{-s}\times H^{-s}$ to the
identified finite projections and this joint tightness proves the
stated Gaussian law on the product space. Its block diagonal covariance
also gives independence of the two limiting coordinates.
\end{proof}

\begin{proof}[Proof of Corollary~\ref{cor:compact-iid}]
The empirical type fluctuation enters through $h_n(\rho\mu)$.  For the fixed
kernel, we have exactly
\[
 h_n(\rho\mu)=-\rho(\mu^{(n)}-\mu),\qquad
 \sqrt n\,h_n(\rho\mu)=-\rho H_n.
\]
The empirical weak convergence also gives $\|\mu^{(n)}-\mu\|_{-s}\to0$
almost surely by Lemma~\ref{lem:compact-analysis}, so
\eqref{eq:compact-derivative} holds in probability.
The approximation \eqref{eq:compact-exploration} and equation
\eqref{eq:exploration-equation} likewise hold under the joint law
by conditioning on the types.  The joint limit in
Lemma~\ref{lem:compact-joint-noise} gives tightness of $H_n$ and $Z_n$,
so the two residuals are $\Op(1)$ after multiplication by $\sqrt n$.
Lemmas~\ref{lem:compact-inverse} and~\ref{lem:linearization} therefore give
\[
 \frac{\mathcal L_n-n\rho\mu}{\sqrt n}
 =J^{-1}(\rho H_n-Z_n)+\op(1).
\]

Let $M_\rho:H^{-s}\to H^{-s}$ denote multiplication by $\rho$, which
is bounded by \eqref{eq:compact-multiplier}.  Apply the bounded linear
map $(\eta,\zeta)\mapsto M_\rho\eta-\zeta$ to the joint Gaussian limit
of Lemma~\ref{lem:compact-joint-noise}.
Its two independent blocks give covariance operator
\begin{align*}
 &M_\rho\{\mathcal V_\mu(1)-\mu\otimes\mu\}
             M_\rho^{\dagger_{-s}}+\mathcal V_\mu(\rho(1-\rho))\\
 &\qquad=\mathcal V_\mu(\rho^2)
             -(\rho\mu)\otimes(\rho\mu)+\mathcal V_\mu(\rho(1-\rho))\\
 &\qquad=\mathcal V_\mu(\rho)-(\rho\mu)\otimes(\rho\mu).
\end{align*}
Here $M_\rho\delta_x=\rho(x)\delta_x$ and
$M_\rho\mu=\rho\mu$, which justify the first equality directly from
\eqref{eq:measure-covariance-notation}.
The resulting operator is positive and trace class, with trace at most
\[
 \operatorname{tr}\mathcal V_\mu(\rho)
 \leq C_s^2\int_M\rho\,\dd\mu<\infty.
\]
The continuous mapping theorem, now applied to $J^{-1}$, and Slutsky's
theorem yield
\[
 \frac{\mathcal L_n-n\rho\mu}{\sqrt n}
 \Rightarrow\Normal\!\left(0,
 J^{-1}\{\mathcal V_\mu(\rho)-(\rho\mu)\otimes(\rho\mu)\}
 (J^{-1})^{\dagger_{-s}}\right),
\]
which is \eqref{eq:compact-iid-limit}.
Using $(J^{-1}\eta)(f)=\eta(Rf)$, its covariance at $f,g\in H^s$ is
\begin{equation}
 \int_M (Rf)(Rg)\rho\,\dd\mu
 -\left(\int_M Rf\,\rho\,\dd\mu\right)
  \left(\int_M Rg\,\rho\,\dd\mu\right).
 \label{eq:compact-iid-covariance}
\end{equation}
\end{proof}

\subsection{Heat kernels and spectral covariances}\label{sec:compact-heat-proof}

For the heat kernel example, $J^*$ is diagonal in the Laplace eigenbasis.
This lets us compute the limiting variance of each mode and verify
trace norm convergence of the covariance series.

\begin{proof}[Verification of Example~\ref{ex:compact-smooth}]
\noindent Fixed point and spectrum.
The heat kernel on a closed manifold is smooth, symmetric, nonnegative,
and preserves constants; see \cite[Chapter~6]{Canzani2013}.  Hence
$\kappa\geq\lambda-\varepsilon>0$, and its row integral is $\lambda$.
The heat kernel expansion \eqref{eq:heat-kernel-expansion}, applied to
\eqref{eq:heat-kernel-model}, gives
\[
 T_\kappa e_0=\lambda e_0,\qquad
 T_\kappa e_j=\varepsilon e^{-u\lambda_j}e_j\quad(j\geq1).
\]
This remains true when $M$ is disconnected: the remaining zero
Laplace eigenfunctions are orthogonal to $e_0=1$.
Thus $\|T_\kappa\|=\lambda>1$.  Since $d\geq1$, there are infinitely
many eigenfunctions, and all the coefficients
$\varepsilon e^{-u\lambda_j}$ are positive; hence $T_\kappa$ has
infinite rank.
To identify the maximal solution of \eqref{eq:compact-fixed-point},
start the fixed point iteration
at $f_0=1$.  If $f_m$ is constant, the row integral gives
$f_{m+1}=1-e^{-\lambda f_m}$, again constant.  The decreasing iterates
converge to the maximal solution.  Since $\lambda>1$, their limit is the
positive scalar solution $\rho=1-e^{-\lambda\rho}$.
For the numbers $r_j$ in \eqref{eq:heat-kernel-modes},
$\lambda=-\log(1-\rho)/\rho$ and $-(1-\rho)\log(1-\rho)<\rho$ give
\[
 (1-\rho)\lambda<1,\qquad
 r_j\geq1-(1-\rho)\lambda>0\quad(j\geq0).
\]

\par\smallskip\noindent Spectral covariance operators.
Because $\rho$ is constant and the kernel is symmetric,
$J^*=I-(1-\rho)T_\kappa$ on $H^s$, and hence
\begin{equation}
 Re_j=(J^*)^{-1}e_j=r_j^{-1}e_j\quad(j\geq0).
 \label{eq:heat-resolvent}
\end{equation}
For iid types, substitute \eqref{eq:heat-resolvent} in
\eqref{eq:compact-iid-covariance}.  Orthonormality and
$\mu(e_j)=\ind_{\{j=0\}}$ give the spectral covariance entries
\[
 \frac{\rho\ind_{\{j=k\}}-
              \rho^2\ind_{\{j=k=0\}}}{r_jr_k}
 =\ind_{\{j=k\}}\sigma_j^2.
\]
For deterministic types with $\beta=\Lambda=0$, the covariance in
Theorem~\ref{thm:compact} similarly gives the entries
\[
 \frac{(1-\rho)\rho}{r_jr_k}\ind_{\{j=k\}}.
\]
These diagonal Gaussian covariances show that every finite collection of
limiting spectral coordinates is independent.  Coordinate $j$ has law
$\Normal_1(0,\sigma_j^2)$ for iid types and
$\Normal_1(0,\rho(1-\rho)/r_j^2)$ for deterministic types.

It remains to check that the covariance series define operators on
$H^{-s}$. Tonelli's theorem and \eqref{eq:compact-delta-tail} give
\[
 \sum_{j\ge0}a_j
 =\int_M\sum_{j\ge0}a_je_j(x)^2\,\mu(\dd x)\leq C_s^2.
\]
For either $d_j=\sigma_j^2$ or $d_j=\rho(1-\rho)/r_j^2$, the
nonnegative coefficients are uniformly bounded because
$r_j\geq1-(1-\rho)\lambda>0$. Since $\|e_j\mu\|_{-s}^2=a_j$,
\[
 \left\|\sum_{j=m}^N d_j(e_j\mu)\otimes(e_j\mu)\right\|_{\mathrm{tr}}
 =\sum_{j=m}^N d_ja_j\longrightarrow0
 \quad\hbox{as }m,N\to\infty.
\]
Thus $\mathcal C_{\mathrm{iid}}$ and $\mathcal C_{\mathrm{det}}$ are
positive trace class operators. Their evaluation covariances are the
entries computed above, so they give the laws
$\Normal(0,\mathcal C_{\mathrm{iid}})$ and
$\Normal(0,\mathcal C_{\mathrm{det}})$ in
\eqref{eq:heat-iid-modes}--\eqref{eq:heat-deterministic-modes}.
The constant mode has the same variance as the constant kernel
$\lambda$: both kernels have row integral $\lambda$, so they have the
same $\rho$ and $r_0=1-(1-\rho)\lambda$. For each nonconstant mode, the ratio
of its variance to that for the constant kernel is
$(1-(1-\rho)\varepsilon e^{-u\lambda_j})^{-2}$, for either sampling scheme;
this amplification factor decreases to one as $\lambda_j\to\infty$.

For \eqref{eq:heat-kernel-perturbation}, boundedness of $e_1$ and
$\kappa\geq\lambda-\varepsilon>0$ ensure that $\kappa_n\geq0$ for all
sufficiently large $n$.  Since $\mu(e_1)=0$, we have
$\mathsf K_\Lambda(\rho\mu)=\rho e_1$ and
$J(e_1\mu)=r_1e_1\mu$.  With $\beta=0$, the mean in
Theorem~\ref{thm:compact} is therefore
$\rho(1-\rho)e_1\mu/r_1$, and its covariance is unchanged.  This gives
\eqref{eq:heat-kernel-perturbed-limit}.

\par\smallskip\noindent The torus grid.
To verify the first order empirical condition for the deterministic grid
stated in the example, take on the flat torus the
regular grid with $n=m^d$ points.
For the complex Fourier characters $\chi_k(x)=e^{2\pi\mathrm i k\cdot x}$,
$k\in\mathbb Z^d$, direct summation over the grid gives
$\mu^{(n)}(\chi_k)=\ind_{\{k\in m\mathbb Z^d\}}$.
Changing from real to complex orthonormal eigenfunctions leaves the
Sobolev norm unchanged within each Laplace eigenspace.  Since
$\mu(\chi_k)=\ind_{\{k=0\}}$, it follows that
\[
 \|\mu^{(n)}-\mu\|_{-s}^2
 =\sum_{k\in\mathbb Z^d\setminus\{0\}}
       (1+4\pi^2m^2|k|^2)^{-s}
 \leq C m^{-2s}.
\]
Thus $\sqrt n\,\|\mu^{(n)}-\mu\|_{-s}\leq C m^{d/2-s}\to0$.
The fixed kernel has $\Lambda=0$, so the deterministic conclusion applies.
\end{proof}

\section*{Acknowledgements}
Part of this work was carried out by Haoxuan Ma during his undergraduate thesis research at the University of Science and Technology of China in 2025–2026. We used AI tools to help organize and refine this manuscript. The authors take full responsibility for the mathematical content of the paper.

\end{document}